\documentclass[11pt]{article}
\usepackage[utf8]{inputenc} 
\usepackage{amssymb, amsmath, amsthm, amsfonts,xcolor, enumerate}
\usepackage{tikz}
\usetikzlibrary{patterns}
\usetikzlibrary{shapes}
\usepackage{fullpage}
\usepackage{hyperref}
\usepackage{ifthen}
\usepackage[style=american]{csquotes}
\usepackage{color}
\usepackage{tikz}  
\usepackage{pgffor} 
\usetikzlibrary{positioning} 
\usepackage{float} 
\usepackage{paralist}

\def\opn#1#2{\def#1{\operatorname{#2}}}
\opn\ext{ex}
\opn\max{max}
\newcommand{\eps}{\varepsilon}
\newcommand{\cG}{\mathcal{G}}
\newcommand{\bN}{\mathbb{N}}
\newcommand{\bR}{\mathbb{R}}
\newcommand{\de}{\delta}
\newcommand{\al}{\alpha}
\newcommand{\ga}{\gamma}
\newcommand{\el}{\ell}
\newcommand{\cM}{\mathcal{M}}

\newcommand{\card}[1]{\left| #1 \right|}

\newtheorem{theorem}{Theorem}[section]
\newtheorem{lemma}[theorem]{Lemma}

\newtheorem{conj}[theorem]{Conjecture}

\newtheorem{claim}[theorem]{Claim}
\newtheorem{prop}[theorem]{Proposition}

\theoremstyle{definition}
\newtheorem{definition}[theorem]{Definition} 
\newtheorem{remark}[theorem]{Remark}

\title{A Density Tur\'an Theorem}
\author{Lothar Narins \thanks{\tt{narins@zedat.fu-berlin.de}}
\and Tuan Tran \thanks{\small \tt{tran@cs.cas.cz}. {\rm Institute of Computer Science, Czech Academy of Sciences, Pod Vodárenskou věží 271/2, 18200 Prague, Czech Republic. Research supported by DFG within the Research Training Group ``Methods for Discrete Structures'' (GRK 1408).}}}
\date{}
\begin{document}
\maketitle

\begin{abstract}
Let $F$ be a graph which contains an edge whose deletion reduces its chromatic number. For such a graph $F,$ a classical result of Simonovits from 1966 shows  that every graph on $n>n_0(F)$ vertices with more than $\frac{\chi(F)-2}{\chi(F)-1}\cdot \frac{n^2}{2}$ edges contains a copy of $F$.
In this paper we derive a similar theorem for multipartite graphs.

For a graph $H$ and an integer $\el \ge v(H)$, let $d_{\el}(H)$ be the minimum real number such that every $\el$-partite graph whose edge density between any two parts is greater than $d_{\el}(H)$ contains a copy of $H$. 
Our main contribution in this paper is to show that $d_{\el}(H)=\frac{\chi(H)-2}{\chi(H)-1}$ for all $\el\ge \el_0(H)$ sufficiently large if and only if $H$ admits a vertex-colouring with $\chi(H)-1$ colours such that all colour classes but one are independent sets, and the exceptional class induces just a matching. When $H$ is a complete graph, this recovers a result of Pfender [Complete subgraphs in multipartite graphs, \emph{Combinatorica} \textbf{32} (2012), 483--495]. We also consider several extensions of Pfender's result.

\end{abstract}

\section{Introduction}
Extremal graph theory has enjoyed tremendous growth in recent decades. One of the central questions from which the theory originated can be described as follows. Given a \emph{forbidden graph} $H$, the Tur\'an problem asks to determine $\ext(n,H)$, the maximum possible number of edges in a graph on $n$ vertices without a copy of $H$. This number is called the \emph{Tur\'an number} of $H$.
Instances of this problem have many connections and applications to other areas. 
In this paper we consider a multipartite version of the problem, suggested by Bollob\'as \cite{Bollobas78}.
Before stating the problem at hand and presenting our contributions, we begin with a brief survey of relevant results.  

\subsection{Background} \label{subsect:turan_problem}

The fundamental Tur\'an theorem of 1941 \cite{Turan41} completely determined the Tur\'an numbers of a clique: the Tur\'an graph $T_{k-1}(n)$, the complete ${(k-1)}$-partite graph on $n$ vertices with parts as equal as possible, has the largest number of edges among all $K_k$-free $n$-vertex graphs. Thus, we have $\ext(n, K_k) = t_{k-1}(n)$, where $t_{k-1}(n)$ is the number of edges in $T_{k-1}(n)$. This theorem generalises a previous result by Mantel \cite{Mantel07} from 1907, which states that $\ext(n, K_3) = \lfloor\frac{n^2}{4}\rfloor$.
 
A large and important class of graphs for which the Tur\'an numbers are well-understood is formed by \emph{colour-critical graphs}, that is, graphs whose chromatic number can be decreased by removing an edge.
Simonovits \cite{Simonovits66} introduced the stability method to show that $\ext(n, H) = t_{k-1}(n)$ for all $n \ge n_0(H)$ sufficiently large, provided $H$ is a colour-critical graph with $\chi(H)=k$; furthermore, $T_{k-1}(n)$ is the unique extremal graph. As the cliques are colour-critical, Simonovits' theorem implies Tur\'an's theorem for large $n$.

For general graphs $H$ we still do not know how to compute the Tur\'an numbers $\ext(n,H)$ exactly; but if we are satisfied with an approximate answer the theory becomes quite simple: it is enough to know the chromatic number of $H$. The important and deep theorem of Erd\H{o}s and Stone \cite{ESt46} together with an observation of Erd\H{o}s and Simonovits \cite{ES66} shows that $\ext(n,H)=\left(\frac{\chi(H)-2}{\chi(H)-1}+o(1)\right)\frac{n^2}{2}$,  
where the $o(1)$ term tends to $0$ as $n$ tends to infinity. In the literature, this result is usually referred as the Erd\H{o}s--Stone--Simonovits theorem.

In the years since these seminal theorems appeared, great efforts have been made to extend them, some of which are discussed in Nikiforov's survey \cite{Nikiforov11}. We are particularly interested in the following two extensions.
  
For every integer $s \ge 2$, let $K_{k-1}(s)$ denote the complete $(k-1)$-partite graph $K_{k-1}(s,\ldots,s)$, and let $K_{k-1}^+(s)$ be the graph obtained from $K_{k-1}(s)$ by adding an edge to the first class. Nikiforov \cite{Nikiforov10} and Erd\H{o}s \cite{Erdos63} (for $k=3$) proved that for all $k \ge 3$ and all sufficiently small $c>0$, every graph of sufficiently large order $n$ with $t_{k-1}(n)+1$ edges contains not only a $K_k$ but a copy of $K_{k-1}^+\bigl(\lfloor c \ln n \rfloor\bigl)$. For fixed $k$, the Erd\H{o}s-R\'enyi random graph $G_{n,p}$ shows that the lower bound $c\ln n$ on the size of the subgraph in this result is optimal up to a constant factor.

Seeking an extension of Tur\'an's theorem, Erd\H{o}s \cite{Erdos69} asked how many $K_k$ sharing a common edge must exist in a graph on $n$ vertices with $t_{k-1}(n)+1$ edges. Bollob\'as and Nikiforov \cite{BN08} sharpened Erd\H{o}s's result \cite{Erdos69} showing that for large enough $n$, every graph of order $n$ with $t_{k-1}(n)+1$ edges has an edge that is contained in $k^{-k-4}n^{k-2}$ copies of $K_k$. This result is 
best possible, up to a $\hbox{poly}(k)$ factor.

In this paper we shall study analogues of these results for multipartite graphs. For a graph $H$ and an integer ${\el \ge v(H)}$, let $d_{\el}(H)$ be the minimum real number such that every $\el$-partite graph $G=(V_1 \cup \ldots \cup V_{\el},E)$ with $d(V_i,V_j):=\frac{e(V_i,V_j)}{\card{V_i}\card{V_j}} > d_{\el}(H)$ for all $i \ne j$ contains a copy of $H$.
The problem of determining the exact value of $d_{\el}(H)$ was suggested by Bollob\'as (see the discussion after the proof of Theorem VI.2.15 in \cite{Bollobas78}). However, it was first studied systematically by Bondy, Shen, Thomass\'e and Thomassen \cite{BSTT06}.
Amongst other things Bondy et.al. showed that for every graph $H$ the sequence $d_{\el}(H)$ decreases to $\frac{\chi(H)-2}{\chi(H)-1}$ as $\el$ tends to infinity. 
To show the lower bound $d_{\el}(H)\ge \frac{\chi(H)-2}{\chi(H)-1}$, they observed that the $\el$-partite graph $G$ obtained from the empty graph on $\{1,\ldots,\el\}$ by splitting each vertex $v$ of $\{1,\ldots,\el\}$ into $\chi(H)-1$ vertices $v_1,v_2,\ldots,v_{\chi(H)-1}$, and joining two vertices $x_i$ and $y_j$ if and only if $x\ne y$ and $i\ne j$, has all edge densities equal to $\frac{\chi(H)-2}{\chi(H)-1}$. Since $G$ is $(\chi(H)-1)$-colourable (with vertex classes $V_i=\{v_i:v\in\{1,\ldots,\el\}\}$ for $1\le i\le \chi(H)-1$), it does not contain a copy of $H$. For the opposite inequality $\lim\limits_{\el \rightarrow \infty}d_{\el}(H)\le \frac{\chi(H)-2}{\chi(H)-1}$, they used the Erd\H{o}s--Stone--Simonovits theorem together with an averaging argument.

When $H=K_3$, the aforementioned result of Bondy et. al. \cite{BSTT06} implies that $d_{\el}(K_3)$ decreases to $\tfrac{1}{2}$ as $\el$ tends to infinity. They also showed that $d_3(K_3)=\frac{-1+\sqrt{5}}{2}\approx 0.61$, $d_4(K_3)>0.51$, and speculated that $d_{\el}(K_3)>\frac{1}{2}$ for all $\el \ge 3$.
Refuting this conjecture, Pfender \cite{Pfender12} proved that $d_{\el}(K_k)=\frac{k-2}{k-1}$ for large enough $\el$. He also described the family $\cG_{\el}^{k}$ of extremal graphs; we shall define this family later in Section \ref{sect:extremal_graphs}.
   
\begin{theorem} [Pfender \cite{Pfender12}]
\label{thm:Pfender}
For every integer $k\ge 3$ there exists a constant $C=C(k)$ such that the following holds for every integer $\el\ge C$. If $G=(V_1\cup\ldots\cup V_{\el},E)$ is an $\el$-partite graph with
\[
d(V_i,V_j) \ge \tfrac{k-2}{k-1} \ \text{for $i\ne j$,}
\]
then either $G$ contains a $K_k$ or $G$ is isomorphic to a graph in $\cG_{\el}^{k}$. In particular, $d_{\el}(K_k)=\tfrac{k-2}{k-1}$ for every $\el \ge C$. 
\end{theorem}

\noindent This theorem can be seen as a multipartite version of the Tur\'an theorem. For an arbitrary graph $H$, Pfender suggested that $d_{\el}(H)$ should be equal to $\frac{\chi(H)-2}{\chi(H)-1}$ for every $\el\ge \el_0(H)$ sufficiently large.

\subsection{Our results}
In this paper we shows that Pfender's suggestion is not quite true. In fact, we characterise those graphs for which the sequence $d_{\el}(H)$ is eventually constant, calling them {\em almost colour-critical}.
 
\begin{figure}[H]
\centering
\includegraphics[width=0.20\linewidth]{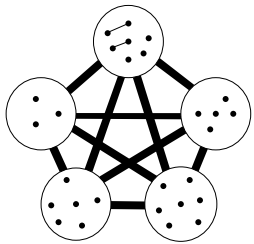}
\caption{An almost colour-critical graph.}
\end{figure}
\vspace{-0.40cm}
\begin{definition}  
A graph $H$ is called almost colour-critical if there exists a map $\phi$ from $V(H)$ to $\{1,2,\ldots,\chi(H)-1\}$ such that
\begin{compactitem}
\item[\rm (i)] The induced subgraph of $H$ on $\phi^{-1}(1)$ has maximum degree at most $1$,
\item[\rm (ii)] For $2 \le i \le \chi(H)-1$, $\phi^{-1}(i)$ is an independent set of $H$.
\end{compactitem}
\end{definition}

\noindent In other words, an almost colour-critical graph $H$ has a vertex-colouring with $\chi(H)-1$ colours that is almost proper: all colour classes but one are independent sets, and the exceptional class induces just a matching (see Figure 1). For example, cliques, or, more generally colour-critical graphs, are almost colour-critical while the complete $k$-partite graphs $K_k(s_1,\ldots,s_k)$ are not for every ${s_1 \ge 1,s_2 \ge 2, \ldots, s_k \ge 2}$.

Our main result shows that almost colour-critical graphs are exactly those for which the sequence $d_{\el}(H)$ is eventually constant.

\begin{theorem} \label{thm:turan}
The following statement holds for every graph $H$. 
\begin{compactitem}
\item[\rm (1)] If $H$ is not almost colour-critical, then $d_{\el}(H) \ge \tfrac{\chi(H)-2}{\chi(H)-1}+ \tfrac{1}{(\chi(H)-1)^2(\el-1)^2}$ for every $\el \ge v(H)$.
\item[\rm (2)] If $H$ is an almost colour-critical graph, then there exists a positive integer $C=C(H)$ so that $d_{\el}(H)=\tfrac{\chi(H)-2}{\chi(H)-1}$ for every $\el >C$.
\end{compactitem}
\end{theorem}

\noindent Note that the estimate in the first statement is tight for $H=K_{1,2}$, and the second statement implies Pfender's result since cliques are almost colour-critical. This result can be viewed as a multipartite version of the Simonovits theorem.  
Since the proof uses the graph removal lemma, the resulting constant $C(H)$ is fairly large. 

The rest of the paper deals with various extensions of Pfender's result. More precisely, we investigate the extensions of Tur\'an's theorem discussed in Section \ref{subsect:turan_problem} for balanced multipartite graphs. An $\el$-partite graph $G$ on non-empty independent sets $V_1,\ldots,V_{\el}$ is {\em balanced} if the vertex classes $V_1,\ldots,V_{\el}$ are of the same size. 

A multipartite version of the extension considered by Nikiforov \cite{Nikiforov10} and Erd\H{o}s \cite{Erdos63} can be stated as follows. 

\begin{theorem} \label{thm:complete_plus-weak}
Let $k$ and $\el$ be integers with $k\ge 3$ and $\el \ge e^{4k^{(k+6)k}}$, and let $G=\left(V_1\cup \ldots \cup V_{\el},E\right)$ be a balanced $\el$-partite graph on $n$ vertices such that
\[
d(V_i,V_j) \ge \tfrac{k-2}{k-1} \quad \text{for $i \ne j$}.
\] 
Then, either $G$ is isomorphic to a graph in $\cG_{\el}^{k}$ or $G$ contains a copy of $K_{k-1}^{+}\bigl(\lfloor c \ln n \rfloor \bigl)$, where $c=k^{-(k+6)k}/2$.
\end{theorem}
\noindent For fixed $k$, the random graph $G_{n,p}$ shows that the lower bound $c\ln n$ on the size of the subgraph in this theorem is tight up to a
constant factor.

The extension of Tur\'an's theorem studied by Bollob\'as and Nikiforov \cite{BN08} has the following multipartite version. 

\begin{theorem} \label{thm:joints}
Let $k$ and $\el$ be integers with $k\ge 3$ and $\el \ge k^{12k}$, and let $G=\left(V_1\cup \ldots \cup V_{\el},E\right)$ be a balanced $\el$-partite graph on $n$ vertices such that
\[
d(V_i,V_j) \ge \tfrac{k-2}{k-1} \quad \text{for $i \ne j$}.
\] 
Then, $G$ either contains a family of $k^{-2k^2}n^{k-2}$ cliques of order $k$ sharing a common edge or is isomorphic to a graph in $\cG_{\el}^{k}$.

\end{theorem}
\noindent With some minor modifications, this result follows from our proof of Theorem \ref{thm:complete_plus-weak}. For the sake of clarity we sketch these modifications after detailing the proof of Theorem \ref{thm:complete_plus-weak}.

\subsection{Organisation}
The remainder of this paper is organised as follows. In Section \ref{sect:preliminaries} we introduce some notation and definitions. In Section \ref{sect:turan_theorem} we extend ideas developed in \cite{Pfender12} to prove Theorem \ref{thm:turan}. A proof of Theorem \ref{thm:complete_plus-weak} is given in Section \ref{sect:complete_plus}. We sketch how to modify the proof of Theorem \ref{thm:complete_plus-weak} to get Theorem \ref{thm:joints} in Section \ref{sect:missing_proofs}, and close with some further remarks and open problems in Section \ref{sect:concluding}.
 
\section{Preliminaries} \label{sect:preliminaries}

\subsection{Notation}
All graphs in this paper are finite, simple and undirected. Given a graph $G$, we denote its vertex and edge sets by $V(H)$ and $E(H)$, and the cardinalities of these two sets by $v(H)$ and $e(H)$, respectively. The minimum degree of $G$ will be denoted by $\de(G)$. 
For a set $U \subseteq V(G)$, we write $G[U]$ for the subgraph of $G$ induced by $U$. The \emph{common neighbourhood} $N(U)$ of $U$ is the set of all vertices of $G$ that are adjacent to every vertex in $U$. Given a vertex $v \in V(G)$, let $\deg(v,U)$ stand for the number of vertices in $U$ adjacent to $v$. For pairwise disjoint vertex sets $W_1,\ldots,W_r \subseteq V(G)$, we write $G[W_1,\ldots,W_r]$ for the $r$-colourable graph which can be obtained from $G[W_1\cup\ldots\cup W_r]$ by deletion of edges in $G[W_i]$ for all $i \le r$. 

Let $G$ be an $\el$-partite graph on non-empty independent sets $V_1,\ldots,V_{\el}$. For $X \subseteq V(G)$ and $i \le \el$, write $X_i=X\cap V_i$. The edge density between $V_i$ and $V_j$ is 
$d_{ij}:=d(V_i,V_j):=\frac{e(V_i,V_j)}{\card{V_i}\card{V_j}}$.

For $r\ge 2$ and $t_1\ge 1,\ldots,t_r\ge 1$, let $K_r(t_1,\ldots,t_r)$ be the complete $r$-partite graph with classes of sizes $t_1,\ldots, t_r$. If $t_1=\ldots=t_r=t$, we simply write $K_r(t)$ instead of $K_r(t_1,\ldots,t_r)$. For $r\ge 2$, $s\ge 1$ and $t_1\ge 2s$, $t_2\ge 1, \ldots,t_r\ge 1$, we denote by $K_r^{+s}(t_1,\ldots,t_r)$ the graph obtained from $K_r(t_1,\ldots,t_r)$ by adding a matching of size $s$ to the first vertex class. If $s=1$, we omit the upper index $s$. In particular, $K_r^{+s}(t)$ is the short form for $K_r^{+s}(t,\ldots,t)$ and $K_r^{+}(t)$ is nothing but $K_r^{+1}(t,\ldots,t)$. 

For $a,b,c \in \bR$, we write $a=b\pm c$ if $b-c\le a\le b+c$. In order to simplify the presentation, we omit floors and ceilings and treat large numbers as integers whenever this does not affect the argument. Unless stated otherwise, all logarithms are base $e$.

The set $\{1,2,\ldots,n\}$ of the first $n$ positive integers is denoted by $[n]$. For $k\in \bN$, we define $\binom{X}{k}:=\{A\subseteq X:\card{A}=k\}$.
We use the symbol $\dot\bigcup$ for union of disjoint sets.

\subsection{Extremal graphs} \label{sect:extremal_graphs}
In this section we shall recall the definition of the family $\mathcal{G}_{\el}^{k}$ of extremal graphs given by Pfender \cite{Pfender12}. For $k\ge 3$ and $\el \ge (k-1)!$, a graph $G$ is in $\bar{\mathcal{G}}_{\el}^{k}$ if it can be constructed as follows. Let $\{\pi_1,\pi_2,\ldots,\pi_{(k-1)!}\}$ be the set of all permutations of $\{1,\ldots,k-1\}$. For $1\le i\le \el$ and $1\le s\le k-1$, pick non-negative integers $n_i^{s}$ such that
\begin{gather*}
n_i^{\pi_i(1)} \ge n_i^{\pi_i(2)}\ge \ldots \ge n_i^{\pi_i(k-1)} \ \text{for} \ 1\le i\le (k-1)!,\\
n_i^{1}=n_i^{2}=\ldots = n_i^{k-1}>0 \ \text{for} \ (k-1)!<i \le \el, \ \textrm{and} \\
\sum_{s} n_i^{s} > 0 \ \text{for} \ 1 \le i \le \el.
\end{gather*}
Vertex and edge sets of $G$ are defined as (see Figure 2)
\begin{align*}
V(G)&=\{(i,s,t): 1 \le i \le \el, 1 \le s \le k-1, 1 \le t \le n_i^{(s)}\},\\
E(G)&=\{(i,s,t)(i',s',t'): i \ne i', s \ne s'\}.
\end{align*}

It is not hard to see that $G$ is an $(k-1)$-colourable $\el$-partite graph with parts $V_i=\{(i,s,t):1\le s\le k-1,1\le t \le n_i^{s}\}$ for $1\le i\le \el$, and colour classes $V^{(s)}=\{(i,s,t):1\le i\le \el,1\le t\le n_i^{s}\}$ for $1\le s\le k-1$. Moreover, if all $n_i^{s}$ are equal, we get $d_{ij}=\frac{k-2}{k-1}$ for every $i \ne j$. Note that other weights $n_i^{(s)}$ can be used to achieve the inequality $d_{ij}\ge \frac{k-2}{k-1}$ for every $i\ne j$.

Let $\mathcal{G}_{\el}^{k}$ be the family of graphs which can be obtained from graphs in $\bar{\mathcal{G}}_{\el}^{k}$ by removal of some edges in $\{(i,s,t)(i',s',t'):1\le i<i'\le (k-1)!\}$.
The following simple observation by Pfender \cite{Pfender12} will be useful for our investigation.

\begin{lemma}
	\label{lem:extremal_graphs}
	Let $k \ge 3$ and $\el \ge (k-1)!$ be integers. If $G=(V_1\cup\ldots\cup V_{\el},E)$ is a $(k-1)$-colourable $\el$-partite graph with $d(V_i,V_j)\ge \tfrac{k-2}{k-1}$ for $i \ne j$, then it is isomorphic to a graph in $\mathcal{G}_{\el}^{k}$.
\end{lemma}

\begin{figure}[H]
\centering
\begin{tikzpicture}[scale=0.70]
\coordinate (1) at (0.75,1); \node [below = 0.17cm of 1] {$n_1^{2}$};
\coordinate (a) at (0.75,3.5); \node [above = 0.17cm of a] {$n_1^{1}$};
\coordinate (2) at (2.25,1); \node [below =0.17cm of 2]{$n_2^{2}$};
\coordinate (b) at (2.25,3.5); \node [above = 0.17cm of b] {$n_2^{1}$};
\coordinate (3) at (3.75,1); \node [below =0.17cm of 3]{$n_3^{2}$};
\coordinate (c) at (3.75,3.5); \node [above = 0.17cm of c] {$n_3^{1}$};
\coordinate (4) at (5.25,1); \node [below =0.17cm of 4]{$\cdots \ n_{\el-1}^{2}$};
\coordinate (d) at (5.25,3.5); \node [above = 0.17cm of d] {$\cdots \ n_{\el-1}^{1}$};
\coordinate (5) at (6.75,1); \node [below =0.17cm of 5]{$\quad n_{\el}^{2}$};
\coordinate (e) at (6.75,3.5); \node [above = 0.17cm of e] {$\quad n_{\el}^{1}$};
\draw (a)--(3)--(b)--(4)--(c)--(5)--(d)--(1)--(e)--(2)--(a);
\draw (a)--(4)--(e)--(3)--(d)--(2)--(c)--(1)--(b)--(5)--(a);
\foreach \x in {(1),(2),(3),(4),(5)}{\draw[fill=blue] \x circle (0.25cm);}
\foreach \y in {(a),(b),(c),(d),(e)}{\draw[fill=red] \y circle(0.25cm);}
\end{tikzpicture}
\caption{A graph in $\bar{\cG}_{\el}^{3}$, all edges between different colours in different parts exists.}
\end{figure}
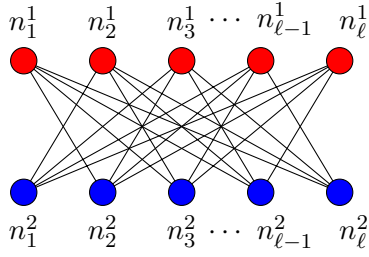

	\subsection{Infracolourable structures}
	The following notation will play a key role in our investigation.
	\begin{definition}
		Given a real number $\eta \ge 0$, and integers $k\ge 3$ and $\el \ge 2$, an $(\eta,k,\el)$-{\em infracolourable structure} is an $\el$-partite graph $G=(V_1\cup\ldots\cup V_{\el},E)$ together with pairs $(D_i^{(s)},Y_i^{(s)})_{s\le k-1,i\le \el}$ satisfying:
		\begin{compactitem}
			\item[\rm (i)] For every $i\le \el$, $V_i=\dot\bigcup_{s\le k-1} Y_i^{(s)}$ and $\card{Y_i^{(1)}}\ge \card{Y_i^{(2)}} \ge \ldots \ge \card{Y_i^{(k-1)}}$;
			\item[\rm (ii)] For every $i\le \el$ and every $s\le k-1$, $D_i^{(s)} \subseteq Y_i^{(s)}$ and $\bigcup_{i\le \el}Y_i^{(s)}\setminus D_i^{(s)}$ is an independent set;
			\item[\rm (iii)] For every $s\le k-1$, each vertex $v\in \bigcup_{i\le \el}D_i^{(s)}$ has at most $\eta \cdot \frac{v(G)}{k-1}$ neighbours in $\bigcup_{i\le \el}Y_i^{(s)}$ and at least $3\eta \cdot \frac{v(G)}{k-1}$ non-neighbours in $\bigcup_{i\le \el}V_i\setminus Y_i^{(s)}$.
		\end{compactitem}
		The graph $G$ is called the \emph{base graph} of the infracolourable structure.
	\end{definition}
	
	Infracolourable structures are useful for us mainly because theirs base graphs break the density conditions in our theorems.
	
	\begin{lemma}\label{lem:infracolourable}
		Let $\eta$ be a positive real number, and let $k\ge 3$ and $\el\ge 2$ be integers. Suppose that an $\el$-partite graph $G=(V_1\cup \ldots\cup V_{\el},E)$ together with a system of pairs $(D_i^{(s)},Y_i^{(s)})_{s\le k-1,i\le \el}$ of vertex sets form an $(\eta,k,\el)$-infracolourable structure. Then
		\[
		e(G)\le \tfrac{k-2}{k-1}\cdot\sum_{i<j}\card{V_i}\card{V_j}.
		\]
		In particular, there exist two different indices $i$ and $j$ such that $d(V_i,V_j)\le \frac{k-2}{k-1}$. Furthermore, the equality occurs if and only if there exists $i_0\in \{0,1,\ldots,\el\}$ such that $D_i^{(s)}=\emptyset$ for all $s$ and all $i$, $\card{Y_i^{(s)}}=\tfrac{1}{k-1}\cdot\card{V_i}$ for all $s$ and all $i\ne i_0$, and $d(Y_i^{(s)},Y_j^{(t)})=1$ for all $s\ne t$ and $i\ne j$. 
	\end{lemma}
	
	\begin{proof}
		It follows from the assumption that
		\begin{gather*}
		e(G)\le\sum_{i<j \atop s \ne t}\card{Y_i^{(s)}}\card{Y_j^{(t)}} + \card{\bigcup_{i,s} D_i^{(s)}}\cdot \left(\eta\cdot\frac{v(G)}{k-1}-\tfrac{1}{2}\cdot 3\eta\cdot\frac{v(G)}{k-1}\right) \\
		\le \sum_{i<j \atop s \ne t}\card{Y_i^{(s)}}\card{Y_j^{(t)}}
		= \sum_{i<j}\card{V_i}\card{V_j}-\sum_{i<j \atop s \le k-1}\card{Y_i^{(s)}}\card{Y_j^{(s)}} \le \tfrac{k-2}{k-1}\cdot \sum_{i<j}\card{V_i}\card{V_j},  
		\end{gather*} 
		where in the last inequality we use Chebyshev's sum inequality.
	\end{proof}
	
	To find an infracolourable structure in host graphs we shall need the following technical lemma. It was implicitly stated in \cite{Pfender12}. We include a proof here for the sake of completeness. 
	
	\begin{lemma} \label{lem:balanced}
		Let $k\ge 3$ and $\el\ge 2$ be integers, and let $\eps$ be a real number with $0<\eps <\tfrac{1}{4}$. Suppose that $G=(V_1\cup\ldots\cup V_{\el},E)$ is an $\el$-partite graph with $d(V_i,V_j)\ge \tfrac{k-2}{k-1}$ for all $i\ne j$. Assume that $X_i^{(s)}$ and $T_i$ be subsets of $V(G)$ for $i\le \el$ and $s\le k-1$ with the following three properties: 
		\begin{compactitem}
			\item[\rm (i)] For every $i\le \el$, $V_i=X_i^{(1)}\dot\cup\ldots\dot\cup X_i^{(k-1)}\dot\cup T_i$;
			\item[\rm (ii)] For every $i\le \el$, $\card{X_i^{(1)}}\ge \ldots \ge \card{X_i^{(k-1)}}$ and $\card{T_i}\le \eps \card{V_i}$;
			\item[\rm (iii)] For every $s\le k-1$, $\bigcup_{i\le \el}X_i^{(s)}$ is an independent set.
		\end{compactitem}
		Then there exists a subset $I_0\in \binom{\bN}{k-1}$ so that $\card{X_i^{(s)}}=\left(\tfrac{1}{k-1}\pm k\sqrt{\eps}\right)\card{V_i}$ for $s\le k-1$ and $i\notin I_0$.
	\end{lemma}
	
	\begin{proof}
		It suffices to show that for each $s\le k-1$ there is at most one index $i\le \ell$ such that 
		$\frac{\card{X_{i}^{(s)}}}{\card{V_i}}>\frac{1}{k-1}+\sqrt{\eps}$.
		Assume to the contrary that $\frac{\card{X_{i}^{(s)}}}{\card{V_i}}\ge \frac{\card{X_{j}^{(s)}}}{\card{V_j}}>\frac{1}{k-1}+\sqrt{\eps}$ for some $s$ and $i\ne j$.
		We first prove that $\frac{\card{X_{i}^{(s)}}}{\card{V_i}}\le 1-\eps$. Otherwise, if $\frac{\card{X_{i}^{(s)}}}{\card{V_i}}>1-\eps$, then 
		\[
		d(V_i,V_j) \le 1-\frac{\card{X_i^{(s)}}}{\card{V_i}}\cdot \frac{\card{X_j^{(s)}}}{\card{V_j}} \le 1-(1-\eps)\left(\tfrac{1}{k-1}+\sqrt{\eps}\right) <\tfrac{k-2}{k-1}
		\] 
		for $k\ge 3$ and $\eps<\tfrac{1}{4}$, as $X_i^{(s)}\cup X_j^{(s)}$ is an independent set by (iii). But this contradicts the density condition that $d(V_i,V_j)\ge \frac{k-2}{k-1}$. 
		
		We shall get a contradiction by proving that $d(V_i,V_j) < \frac{k-2}{k-1}$. Indeed, we can infer from Chebyschev's sum inequality that 
		\begin{align*}
		d(V_i,V_j) & \overset{(iii)}{\le} 1 -\frac{1}{\card{V_i}\card{V_j}}\cdot \sum_{t}\card{X_i^{(t)}}\card{X_j^{(t)}} \\
		& \le     1-\frac{\card{X_i^{(s)}}\card{X_j^{(s)}}}{\card{V_i}\card{V_j}}-\frac{1}{(k-2)\card{V_i}\card{V_j}}
		\cdot\left(\card{V_i}-\card{T_i}-\card{X_i^{(s)}}\right)\left(\card{V_j}-\card{T_j}-\card{X_j^{(s)}}\right)\\
		&= 1-x_ix_j-\tfrac{1}{k-2}\left(1-t_i-x_i\right)\left(1-t_j-x_j\right),
		\end{align*}
		where $x_i=\frac{\card{X_{i}^{(s)}}}{\card{V_i}}$, $x_j=\frac{\card{X_{j}^{(s)}}}{\card{V_j}}$, $t_i=\frac{\card{T_i}}{\card{V_i}}$ and $t_j=\frac{\card{T_j}}{\card{V_j}}$. Since both $x_i$ and $x_j$ are bounded from below by $\frac{1}{k-1}$, the expression $f(x_i,x_j,t_i,t_j):=1-x_ix_j-\tfrac{1}{k-2}\left(1-t_i-x_i\right)\left(1-t_j-x_j\right)$ is decreasing with respect to both $x_i$ and $x_j$. Therefore, the density $d(V_i,V_j)$ is bounded from above by           
		\begin{align*}
		f(x_i,x_j,t_i,t_j) &\le f\left(\tfrac{1}{k-1}+\sqrt{\eps},\tfrac{1}{k-1}+\sqrt{\eps},t_i,t_j\right)\le f\left(\tfrac{1}{k-1}+\sqrt{\eps},\tfrac{1}{k-1}+\sqrt{\eps},\eps,\eps\right)  < \tfrac{k-2}{k-1},
		\end{align*}
		where the second inequality follows from the assumption that $t_i,t_j \in [0,\eps]$. However, this contradicts the assumption that $d(V_i,V_j) \ge \frac{k-2}{k-1}$.
	\end{proof}

	\section{Proof of Theorem \ref{thm:turan}} \label{sect:turan_theorem}
	In this section we will prove Therem \ref{thm:turan}. We begin with a proof of the first assertion.
	
	\begin{proof}[Proof of Theorem \ref{thm:turan}(1)] 
We prove by contradiction. Assume that $d_{\el}(H)<\frac{\chi(H)-2}{\chi(H)-1}+\frac{1}{(\chi(H)-1)^2(\el-1)^2}$. Let $r=\chi(H)-1$, and let $V_1,\ldots,V_{\el}$ be $\el$ disjoint sets of size $(\el-1)r$. For $i \le \el$, we partition $V_i$ into $r$ subsets $V_i^{(1)},\ldots,V_i^{(r)}$ of size $(\el-1)$ each. We form a complete bipartite graph between $V_i^{(s)}$ and $V_j^{(t)}$ for $i<j$ and $s \ne t$. We then create a perfect matching in $V_1^{(1)} \cup \ldots \cup V_{\el}^{(1)}$ such that there is exactly one edge between $V_i^{(1)}$ and $V_j^{(1)}$ for every $i \ne j$. 
		The resulting graph $G$ satisfies
		\[
		d(V_i,V_{j})=\frac{\chi(H)-2}{\chi(H)-1}+\frac{1}{(\chi(H)-1)^2(\el-1)^2} >d_{\el}(H)\quad \text{for $i \ne j$}.
		\]
		Thus, by the definition of $d_{\el}(H)$, $G$ must contain a copy of $H$. From the construction of $G$, we can see that $H$ is an almost colour-critical graph. This finishes our proof of Theorem \ref{thm:turan}(1).
	\end{proof}
	
	\begin{remark}\label{rmk:tighness}
		The estimate in Theorem \ref{thm:turan}(1) is tight for $K_{1,2}$, that is $d_{\el}(K_{1,2})=\frac{1}{(\el-1)^2}$ for $\el \ge 3$. Indeed, let $G=(V_1\cup\ldots\cup V_{\el},E)$ be an $\el$-partite graph with $d(V_i,V_j)>\frac{1}{(\el-1)^2}$ for every $i\ne j$. We wish to show that $G$ contains a copy of $K_{1,2}$. Suppose to the contrary that $G$ is $K_{1,2}$-free. For $i\ne j$, we write $V_{i,j}$ for the set of vertices in $V_i$ with at least one neighbour in $V_j$. 
		Since $G$ is $K_{1,2}$-free, we see that
		\begin{itemize}
			\item[\rm(i)] the edges between $V_i$ and $V_j$ form a perfect matching between $V_{i,j}$ and $V_{j,i}$ for every $i\ne j$;
			\item[\rm (ii)] $V_{i,j}$ and $V_{i,j'}$ are disjoint for all distinct indices $i,j$ and $j'$.
		\end{itemize} 
		Notice that $V_{i,j}$ is non-empty for every $i\ne j$ as $d(V_i,V_j)>0$. Combining this with property (ii), we conclude that 
		\begin{equation}\label{eq:V-lower}
		\card{V_i} \ge \sum_{j\in [\el]\setminus \{i\}}\card{V_{i,j}} \ge \el-1 \ \text{for $i\le \el$}.
		\end{equation}
		Hence
		\[
		\sum_{1\le i<j \le \el}\left(\frac{\card{V_{i,j}}}{\card{V_i}}+\frac{\card{V_{j,i}}}{\card{V_j}}\right)=\sum_{1\le i\le \el}\left(\sum_{j'\ne i}\frac{\card{V_{i,j'}}}{\card{V_i}}\right) \le \el.
		\]
		Consequently, there exist $1\le i<j\le \el$ with $\frac{\card{V_{i,j}}}{\card{V_i}}+\frac{\card{V_{j,i}}}{\card{V_j}} \le \frac{\el}{\binom{\el}{2}}=\frac{2}{\el-1}$. By appealing to the AM-GM inequality, we thus get $\sqrt{\card{V_{i,j}}\card{V_{j,i}}} \le \frac{1}{\el-1}\cdot \sqrt{\card{V_i}\card{V_j}}$. This forces
		\[
		d(V_i,V_j)\overset{(i)}{=} \frac{\card{V_{i,j}}}{\card{V_i}\card{V_j}}\overset{(i)}{=} \frac{\sqrt{\card{V_{i,j}}\card{V_{j,i}}}}{\card{V_i}\card{V_j}} \le \frac{1}{(\el-1)\sqrt{\card{V_i}\card{V_j}}} \overset{\eqref{eq:V-lower}}{\le} \frac{1}{(\el-1)^2},
		\]
		contradicting the assumption that $d(V_i,V_j)>\frac{1}{(\el-1)^2}$.
	\end{remark}
	
	To handle the second statement of Theorem \ref{thm:turan}, we shall prove a stronger result.  
	\begin{theorem} \label{thm:key}
		Let $H$ be an almost colour-critical graph. Then, there exists a constant $C=C(H)$ such that for every integer $\el>C$, every $\el$-partite graph $G=(V_1\cup \ldots \cup V_{\el},E)$ with
		\[
		d(V_i,V_j) > \frac{\chi(H)-2}{\chi(H)-1} \quad \text{for $i \ne j$}
		\]
		contains a copy of $H$ whose vertices are in different parts of $G$.
	\end{theorem}
	
	\begin{remark} \label{rmk:universal_graph}
		Suppose that $H$ is almost colour-critical. Let $k=\chi(H)$ and $q=v(H)$. From the definition of almost colour-critical graphs, $H$ is a subgraph $K^{+q}_{k-1}(2q)$. Moreover, it is easy to see that $\chi(K^{+q}_{k-1}(2q))=k=\chi(H)$ and $K^{+q}_{k-1}(2q)$ is almost colour-critical. Therefore, if Theorem \ref{thm:key} holds for $K^{+q}_{k-1}(2q)$, it will hold for $H$ as well.
	\end{remark}
	
	The main idea of the proof of Theorem \ref{thm:key} is as follows. Let $G=(V_1\cup\ldots\cup V_{\el},E)$ be a counterexample. We first apply a stability result (Lemma \ref{lem:induced_subgraph}) to obtain an induced $(\chi(H)-1)$-colourable subgraph of $G$ which almost spans $V(G)$. Using embedding results (Lemmas \ref{lem:common_neighbours} and \ref{lem:embedding_complete}) we can then show that there exists a subset $I\subseteq [\el]$ such that $G[\bigcup_{i\in I}V_i]$ is the base graph of an $(\eta,k,\card{I})$-infracolourable structure. But according to Lemma \ref{lem:infracolourable}, this forces $d(V_i,V_j)\le \frac{k-2}{k-1}$ for some $i,j \in I$, violating the density condition. 
	
	Our first step in the proof of Theorem \ref{thm:key} will be to show that a counterexample $G$ must contain an induced $(\chi(H)-1)$-colourable subgraph which almost spans $V(G)$. For that we shall need the following stability result.
	
	\begin{lemma} \label{lem:induced_subgraph}
		Given integers $k\ge 3$ and $q\ge 1$ and a real number $0<\eps<\frac{1}{8k^2q}$, there exists a constant $C=C(k,q,\eps)$ such that the following holds for $\el \ge C$. Let $G=(V_1\cup\ldots\cup V_{\el},E)$ be a balanced $\el$-partite graph on $n$ vertices with $d(V_i,V_j)\ge \frac{k-2}{k-1}$ for all $i\ne j$. Suppose $G$ contains no copy of $K_{k-1}^{+q}(2q)$ whose vertices lie in different parts of $G$. Then, $G$ contains an induced $(k-1)$-colourable subgraph $F$ whose vertex classes $X^{(1)}, \ldots, X^{(k-1)}$ satisfy the following properties   
		\begin{compactitem}
			\item[\rm (i)] For $s \le k-1$, $\card{X^{(s)}}= \left(\tfrac{1}{k-1}\pm\eps\right)n$;
			\item[\rm (ii)] For $s \le k-1$ and $v\in \bigcup_{t\ne s}X^{(t)}$, $\deg(v,X^{(s)})\ge \card{X^{(s)}}-\eps n$. 
		\end{compactitem}
	\end{lemma}
	
	To prove Lemma \ref{lem:induced_subgraph} we require the following result whose proof can be found in Section \ref{sect:missing_proofs}.
	
	\begin{prop} \label{prop:Erdos-Stone_stability}
		For every graph $H$ and every $\eps>0$, there exist positive constants $\ga=\ga(H,\eps)$ and $C=C(H,\eps)$ such that the following holds for $n \ge C$. Suppose that $G$ is an $n$-vertex graph with
		${e(G) \ge \left(\frac{\chi(H)-2}{\chi(H)-1}-\ga \right)\binom{n}{2}}$ containing at most $\ga n^{v(H)}$ copies of $H$. Then, $G$ contains a ${(\chi(H)-1)}$-colourable subgraph of order at least $(1-\eps)n$ and minimum degree at least $\left(\frac{\chi(H)-2}{\chi(H)-1}-\eps\right)n$.
	\end{prop}
	
	Another tool that will be used in the proof of Lemma \ref{lem:induced_subgraph} and Theorem \ref{thm:key} is an embedding result. Before stating it, we shall introduce the necessary terminology. Let $G[W^{(1)},\ldots,W^{(r)}]$ be an $r$-colourable graph such that $W^{(s)}=\dot\bigcup_{i\ge 1} W^{(s)}_i$ for every $s\le r$. We call an embedding $f:K_r(a_1,\ldots,a_r) \rightarrow G$ \emph{good} if the $s$th vertex class of $K_r(a_1,\ldots,a_r)$ is mapped to $W^{(s)}$ for every $s \le r$, and for each index $i$ there is at most one vertex $v \in K_r(a_1,\ldots,a_r)$ with $f(v) \in \bigcup_{s\le r}W^{(s)}_i$. 
	
	\begin{lemma} \label{lem:embedding_complete}
		Suppose that $r\ge 2$ and $q\ge 1$ are integers, and let $G[W^{(1)},\ldots,W^{(r)}]$ be an $r$-colourable graph which satisfies the following properties
		\begin{compactitem}
			\item[{\rm (i)}] For $s\le r$, $W^{(s)}=\dot\bigcup_{i}W^{(s)}_i$ and $\card{W^{(s)}_i} < \frac{1}{2rq}\cdot\card{W^{(s)}}$ for all $i$, 
			\item[\rm (ii)] For $s\le r$ and $v \in \bigcup_{t \ne s}W^{(t
				)}$, $\deg(v,W^{(s)}) > (1-\frac{1}{2rq})\cdot\card{W^{(s)}}$.
		\end{compactitem}
		Then, for every $r$-tuple of integers $a_1,\ldots,a_r \in [0,q]$, every good embedding from $K_r(a_1,\ldots,a_r)$ to $G$ can be extended to a good embedding from $K_r(q)$ to $G$.
	\end{lemma}
	\begin{proof} 
		Suppose $f$ is a good embedding from $K_r(a_1,\ldots,a_r)$ to $G$. To prove the lemma, it suffices to show that $f$ can be extended to a good embedding $g$ from $K_r(a_1,\ldots,a_s+1,\ldots,a_r)$ to $G$ whenever $a_s\le q-1$. Let $v$ be the vertex of $K_r(a_1,\ldots,a_s+1,\ldots,a_r)$ which is not in $K_r(a_1,\ldots,a_r)$, and let $X$ denote the set of vertices of $K_r(a_1,\ldots,a_r)$ which are not in the $s$th vertex class. By property (ii), we see that each vertex of $X$ has at most $\frac{1}{2rq}\cdot\card{W^{(s)}}$ non-neighbours in $W^{(s)}$, and thus $\card{N(X)\cap W^{(s)}} \ge \card{W^{(s)}}-\card{X}\cdot\frac{\card{W^{(s)}}}{2rq} \ge \tfrac{1}{2}\card{W^{(s)}}$. Note that, by property (i), each vertex of $X$ can forbid at most $\frac{1}{2rq}\cdot\card{W^{(s)}}$ vertices of $W^{(s)}$ from being the image of $v$. Therefore, the number of possible images of $v$ under $g$ is at least
		$\card{N(X)\cap W^{(s)}}-\card{X}\cdot\frac{\card{W^{(s)}}}{2rq}
		\ge \tfrac{1}{2}\card{W^{(s)}}-\card{X}\cdot\frac{\card{W^{(s)}}}{2rq}>0$,
		where in the last inequality we use the inequality $\card{W^{(s)}}>0$ which is implied by property (i).
	\end{proof}
	
	\begin{proof}[Proof of Lemma \ref{lem:induced_subgraph}]
		We denote $H=K_{k-1}^{+q}(2q)$, and let 
		\begin{equation*}\label{eq:induced_1}
		\ga=\ga_{\ref{prop:Erdos-Stone_stability}}\left(H,\tfrac{\eps}{2k}\right), \ C=\max\biggl\{2k^2q^2\ga^{-1},8(k-1)^2q,4(k-1)q\eps^{-1}, C_{\ref{prop:Erdos-Stone_stability}}(H,\tfrac{\eps}{2k})\biggl\}.
		\end{equation*}
		Because $G=(V_1\cup\ldots\cup V_{\el},E)$ is a balanced $\el$-partite graph on $n$ vertices, we must have
		\begin{equation}
		\label{eq:induced_2}
		\card{V_1}=\card{V_2}=\ldots=\card{V_{\el}}=\frac{n}{\el}:=m.
		\end{equation}
		
		In the first step, we shall use Proposition \ref{prop:Erdos-Stone_stability} to show that $G$ contains an almost spanning $(k-1)$-colourable subgraph.
		Indeed, by the choice of $C$ we see that $n\ge \el \ge C\ge C_{\ref{prop:Erdos-Stone_stability}}(H,\tfrac{\eps}{2k})$.
		Moreover, since $G$ contains no copy of $H$ whose vertices lie in different parts of $G$, the number of copies of $H$ in $G$ is at most 
		\[
		\binom{v(H)}{2}\el m^2n^{v(H)-2} <\tfrac{2k^2q^2}{\el} \cdot (\el m)^2n^{v(H)-2} \le \ga n^{v(H)}, 
		\]
		since $n=\el m$ and $\el \ge C\ge 2k^2q^2\ga^{-1}$.
		Also, by the density condition
		\[
		e(G) \ge \binom{\el}{2}\tfrac{k-2}{k-1}m^2 \overset{\eqref{eq:induced_2}}{\ge} \left(\tfrac{k-2}{k-1}-\tfrac{1}{\el} \right)\binom{n}{2} \ge \left(\tfrac{k-2}{k-1}-\ga\right)\binom{n}{2},
		\]
		assuming $\el \ge C \ge 2k^2q^2\ga^{-1}$.
		Therefore, we can derive from Proposition \ref{prop:Erdos-Stone_stability} that $G$ contains a $(k-1)$-colourable subgraph $F'$ with
		\begin{equation}
		\label{eq:induced_3}
		v(F')\ge(1-\tfrac{\eps}{2k})n \ \textrm{and} \ \de(F')\ge\left(\tfrac{k-2}{k-1}-\tfrac{\eps}{2k}\right)n.
		\end{equation} 
		If $W^{(1)},\ldots, W^{(k-1)}$ are vertex classes of $F'$, then \eqref{eq:induced_3} implies that
		\begin{equation}
		\label{eq:induced_4}
		\left(\tfrac{1}{k-1}-\tfrac{\eps}{2}\right)n\le \card{W^{(s)}} \le \left(\tfrac{1}{k-1}+\tfrac{\eps}{2k}\right)n \quad \textrm{for $s \le k-1$}.
		\end{equation}
		
		In the second step, we shall prove that the induced subgraph $G[V(F')]$ of $G$ does not contain  a large monochromatic matching whose vertices are in different parts of $G$. Indeed, for $s\le k-1$, let $\cM_{(s)}$ denote a maximum matching in $G[W^{(s)}]$ whose vertices are in different parts of $G$, and let $K$ be a subset of $[\el]$ containing all indices $i$ such that $\bigcup_{s\le k-1}\cM_{(s)}$ has a vertex in $V_i$. The size of $K$ will be bounded from above in terms of $k$ and $q$.
		
		\begin{claim}
		\label{eq:induced_5}
		$\card{K} < 2(k-1)q$. 
		\end{claim}
\begin{proof}
We prove the claim by contradiction. Suppose that for some $s \le k-1$, $\cM_{(s)}$ contains a matching of size $q$, say $\{x_1x_2,\ldots,x_{2q-1}x_{2q}\}$, .
We wish to show that the following two properties holds:
		\begin{compactitem}
			\item[\rm (i)] For $t\le k-1$ and $i\le \el$, $W^{(t)}=W_1^{(t)} \dot\cup \ldots \dot\cup W_{\el}^{(t)}$ and $\card{W_i^{(t)}}<\frac{1}{4(k-1)q}\cdot\card{W^{(t)}}$;
			\item[\rm (ii)] For $t\le k-1$ and $v \in V(F')\setminus W^{(t)}$, $\deg_{F'}(v,W^{(t)}) >\left(1-\frac{1}{4(k-1)q}\right)\cdot\card{W^{(t)}}$.
		\end{compactitem}
        Property (i) follows from the estimate
		\[
		\card{W_i^{(t)}} \le \card{V_i}=\frac{n}{\el}<\frac{1}{4(k-1)q}\cdot\left(\frac{1}{k-1}-\tfrac{\eps}{2}\right)n\overset{\eqref{eq:induced_4}}{<} \frac{1}{4(k-1)q}\cdot\card{W^{(t)}}
		\]
		for $\el \ge C \ge 8(k-1)^2q$ and $\eps<\frac{1}{8k^2q}$. To prove (ii), assume that $v\in W^{(s)}$ for some $s\ne t$. Because $W^{(s)}$ is an independent set in $F'$, one has $\card{W^{(t)}}-d_{F'}(v,W^{(t)}) \le v(F')-\card{W^{(s)}}-\deg_{F'}(v)$. Hence by appealing to \eqref{eq:induced_3} and \eqref{eq:induced_4}, we get
		\begin{align*}
		\card{W^{(t)}}-d_{F'}(v,W^{(t)}) &\le n-\left(\frac{1}{k-1}-\tfrac{\eps}{2}\right)n-\left(\frac{k-2}{k-1}-\tfrac{\eps}{2k}\right)n\\
		&\le \eps n<\frac{1}{4(k-1)q}\cdot\left(\frac{1}{k-1}-\tfrac{\eps}{2}\right)n \le \frac{1}{4(k-1)q}\cdot\card{W^{(t)}}
		\end{align*}
		for $\eps<\frac{1}{8k^2q}$. This finishes our verification of (i) and (ii).
		
		Finally, properties (i) and (ii) ensure that we can apply Lemma \ref{lem:embedding_complete} with $r_{\ref{lem:embedding_complete}}=k-1$ and $q_{\ref{lem:embedding_complete}}=2q$ to $G[W^{(1)},\ldots,W^{(r)}]$ to find a copy of $K_{k-1}(2q)$ whose $s$th vertex class is $\{x_1,\ldots,x_{2q}\}$ and vertices lie in different parts of $G$. Since $\{x_1,x_2,\ldots,x_{2q-1}x_{2q}\}$ is a matching in $G$, the graph $G$ contains a desired copy of $H$, which contradicts our hypothesis.
\end{proof}		
		
		To finish the proof, we shall show that $G$ contains an induced subgraph $F$ with the desired properties. For this purpose, we let $X^{(s)}=W^{(s)}\setminus \bigcup_{i \in K}V_i$ for $s\le k-1$. The maximality of $\cM_{(s)}$ implies that $X^{(s)}$ is an independent set in $G$. So the induced subgraph $F=G[X^{(1)} \cup \ldots \cup X^{(k-1)}]$ is $(k-1)$-colourable. What is left is to prove that $F$ has the desired properties. Since $\eps<\frac{1}{8k^2q}$ and $\el\ge C\ge 4(k-1)q\eps^{-1}$, we find that
		\begin{gather*}
		v(F)\ge v(F')-\card{\bigcup_{i\in K}V_i}\overset{\eqref{eq:induced_3},\text{Claim \ref{eq:induced_5}}}{\ge}\left(1-\tfrac{\eps}{2k}\right)n-2(k-1)q\cdot \frac{n}{\el}> (1-\eps)n, \\
		\de(F)\ge \de(F')-\card{\bigcup_{i\in K}V_i}\overset{\eqref{eq:induced_3},\text{Claim \ref{eq:induced_5}}}{\ge}\left(\tfrac{k-2}{k-1}-\tfrac{\eps}{2k}\right)n-2(k-1)q\cdot\frac{n}{\el}>\left(\tfrac{k-2}{k-1}-\tfrac{\eps}{2}\right)n.
		\end{gather*}
		Moreover, by \eqref{eq:induced_4} we see that
		$\card{X^{(s)}} \le \card{W^{(s)}} \le \left(\tfrac{1}{k-1}+\tfrac{\eps}{2k}\right)n$ for $s\le k-1$,
		and hence $\left(\frac{1}{k-1}-\tfrac{\eps}{2}\right)n \le \card{X^{(s)}}\le \left(\frac{1}{k-1}+\tfrac{\eps}{2k}\right)n$ for $s\le k-1$.
		Therefore, for $s\le k-1$ and $v\in \bigcup_{t\ne s}X^{(t)}$, there are at most 
		$n-\card{X^{(s)}}-d_{F}(v) \le n-\left(\tfrac{1}{k-1}-\tfrac{\eps}{2}\right)n-\left(\tfrac{k-2}{k-1}-\tfrac{\eps}{2}\right)n=\eps n$ missing edges in $F$ between $v$ and $X^{(s)}$. This completes our proof of Lemma \ref{lem:induced_subgraph}.
	\end{proof}
	
	We also need the following elementary lemma. It is probably well-known, but we could not find a reference. For completeness we include its proof in Section \ref{sect:missing_proofs}. 
	
	\begin{lemma} \label{lem:common_neighbours}
		Given integers $r\ge 1$ and $q\ge 2$ and a real number $d\in (0,1)$, there exist an integer $D=D(r,q,d)$ and a positive $\rho=\rho(r,q,d)$ so that the following holds. Suppose that $G$ is an $(r+1)$-colourable graph with vertex classes $U,W_{(1)},\ldots,W_{(r)}$. If $\card{U} \ge D$ and $\deg(u,W_{(s)})\ge d\card{W_{(s)}}$ for all $u \in U$ and $s\le r$, then there is a subset $A \in \binom{U}{q}$ with $\card{N(A) \cap W_{(s)}} \ge \rho\card{W_{(s)}}$ for $s\le r$.
	\end{lemma}
	
	To find an infracolourable structure in $G$ we shall make use of a consequence of Lemmas \ref{lem:embedding_complete} and \ref{lem:common_neighbours}.
	
	\begin{lemma}\label{lem:bound_I}
		Given integers $k\ge 3$ and $q \ge 1$ and a real number $\eta \in (0,1)$, there exist integers $C=C(k,q,\eta)$ and $D=D(k,q,\eta)$ and a positive $\de=\de(k,q,\eta)$ such that the following holds for $\el \ge C$ and $\eps \in (0,\de)$. Suppose that $G=(V_1\cup\ldots\cup V_{\el},E)$ is a balanced $\el$-partite graph containing no copy of $K_{k}(2q)$ in $G$ whose vertices are in different parts of $G$. Assume $(X^{(s)}_i)_{s\le k-1, i\le \el}$ are vertex sets satisfying:
		\begin{compactitem}
			\item[\rm (i)] For $i\le \el$, $X_i^{(1)},\ldots, X_i^{(k-1)}$ are disjoint subsets of $V_i$,
			\item[\rm (ii)] For $i\le \el$ and $s\le k-1$, $\card{X^{(s)}_i}= \left(\frac{1}{k-1}\pm\eps\right)\card{V_i}$,
			\item[\rm (iii)] For every $s \le k-1$ and $v\in \bigcup_{i\le \el, t\ne s}X_i^{(t)}$, $\deg(v,\bigcup_{i\le\el}X_i^{(s)}) \ge \card{\bigcup_{i\le \el}X_i^{(s)}}-\eps \cdot v(G)$. 
		\end{compactitem}
		Let $I$ be the subset of $[\el]$ consisting of all indices $i\in [\el]$ such that $V_i$ contains a vertex $v$ with $ \deg(v,\bigcup_{j\le \el}X_j^{(s)}) \ge \eta \cdot v(G)$ for $s \le k-1$. Then $\card{I} \le D$.
	\end{lemma}
	\begin{proof}
		Let  
		$D=D_{\ref{lem:common_neighbours}}\left(k-1,2q,\tfrac{k\eta}{4}\right)$, $C=\max\biggl\{4kD,2\eta^{-1}D,\frac{9(k-1)kq}{\rho}\biggl\}$ and 
		$\de=\min\biggl\{\frac{1}{4k},\frac{\rho}{8(k-1)kq}\biggl\}$, where $\rho=\rho_{\ref{lem:common_neighbours}}\left(k-1,2q,\tfrac{k\eta}{4}\right)$.
		We shall prove the lemma by contradiction. Assume that $\card{I} \ge D$. Let $J$ be an arbitrary subset of $I$ of size $D$. By the definition of $I$, for each index $j\in J$ we can find a vertex $v_j\in V_j$ such that $ \deg(v_j,\bigcup_{i\le \el}X_i^{(s)}) \ge \eta \cdot v(G)$ for $s \le k-1$. Let $U=\{v_j:j\in J\}$. 
		
		For simplicity of notation, let $X^{(s)}:=\bigcup_{i\le \el}X_i^{(s)}$ and $W^{(s)}:=\bigcup_{i \in [\el]\setminus J}X_i^{(s)}$ for $s\le k-1$. Then, property (i) implies that $W^{(1)},\ldots,W^{(k-1)}$ are disjoint subsets of $V(G)$. By (i) and (ii), we find that
		\begin{equation} \label{eq:bound_I_1}
		\card{W^{(s)}}\ge  \left(\frac{1}{k-1}-\eps-\frac{D}{\el}\right)\cdot v(G)\ge \tfrac{v(G)}{2k} 
		\end{equation}
		for $\eps \le \de \le \frac{1}{4k}$ and $\el\ge C\ge 4kD$.
		Also, (i) and (ii) force $\card{W^{(s)}}\le  \left(\frac{1}{k-1}+\eps \right)v(G) \le \frac{2v(G)}{k}$, since $\eps \le \de \le \frac{1}{4k}$.
		Combining these two inequalities, we conclude that 
		\[
		\deg(v,W^{(s)}) \ge \deg(v,X^{(s)})-\card{\bigcup_{j\in J}V_j} \ge \eta \cdot v(G)-D\cdot\frac{v(G)}{\el} \ge \tfrac{\eta}{2}\cdot v(G) \ge \tfrac{k\eta}{4}\cdot \card{W^{(s)}}
		\]
		for $v \in U$ and $s \le k-1$, as $\el \ge 2\eta^{-1}D$. Furthermore, $\card{U}=D= D_{\ref{lem:common_neighbours}}\left(k-1,2q,\tfrac{k\eta}{4}\right)$, by the definition of $D$.  
		By applying Lemma \ref{lem:common_neighbours} to $G[U, W^{(1)}, \ldots, W^{(k-1)}]$ with $r_{\ref{lem:common_neighbours}}=k-1$, $q_{\ref{lem:common_neighbours}}=2q$ and $d_{\ref{lem:common_neighbours}}=\frac{k\eta}{4}$, we thus obtain a subset $A \in \binom{U}{2q}$ with 
		\begin{equation}\label{eq:bound_I_2}
		\card{N(A)\cap W^{(s)}}\ge \rho \card{W^{(s)}}\quad\text{for $s\le k-1$}. 
		\end{equation}

		In the rest of the proof we shall use Lemma \ref{lem:embedding_complete} to show that $G[N(A)\cap W^{(1)},\ldots,N(A)\cap W^{(k-1)}]$ contains a copy of $K_{k-1}(2q)$ whose vertices are in different parts of $G$. Since this copy lies in $N(A)$, together with vertices of $A$ it forms a copy of $K_k(2q)$ whose vertices belong to different parts of $G$, contradicting the assumption. It remains to verify the assumptions of Lemma \ref{lem:embedding_complete}.
		Indeed, for $s \le k-1$, $N(A)\cap W^{(s)}$ does admit the partition    
		\begin{equation}\label{eq:bound_I_3}
		N(A)\cap W^{(s)}=\dot\bigcup_{j \notin J}\bigl(N(A) \cap X_j^{(s)}\bigl).
		\end{equation}
		Moreover, since $N(A)\cap W^{(s)} \subseteq X^{(s)}$ for $s\le k-1$, we must have, for $s \le k-1$ and $v\in \bigcup_{t\ne s}\bigl(N(A)\cap W^{(t)}\bigl)$, 
		\begin{align*}
		\card{N(A)\cap W^{(s)}}-\deg(v,N(A)\cap W^{(s)})&\le \card{\bigcup_{i\le \el}X_i^{(s)}}-\deg(v,\bigcup_{i\le\el}X_i^{(s)})\\
		&\overset{(iii)}{\le}\eps \cdot v(G) \le \tfrac{1}{4(k-1)q} \cdot \rho \cdot \frac{v(G)}{2k} \overset{\eqref{eq:bound_I_1}, \eqref{eq:bound_I_2}}{\le}\tfrac{1}{4(k-1)q}\cdot \card{N(A)\cap W^{(s)}},
		\end{align*}
		assuming $\eps \le \de \le \frac{\rho}{8(k-1)kq}$. It can be rewritten as
		\begin{equation}\label{eq:bound_I_4}
		\deg(v,N(A)\cap W^{(s)}) \ge \left(1-\tfrac{1}{4(k-1)q}\right)\card{N(A)\cap W^{(s)}} \ \text{for $s\le k-1$ and $v\notin \bigcup_{t\ne s}(N(A)\cap W^{(t)}$}.
		\end{equation} 
		Also, for every $j \notin J$ and $s \le k-1$, we have
		\begin{equation}\label{eq:bound_I_5}
		\card{N(A) \cap X_j^{(s)}} \le \card{V_j}=\tfrac{v(G)}{\el}<\tfrac{1}{4(k-1)q} \cdot \rho \cdot \frac{v(G)}{2k} \overset{\eqref{eq:bound_I_1}, \eqref{eq:bound_I_2}}{\le}\tfrac{1}{4(k-1)q}\cdot \card{N(A)\cap W^{(s)}}
		\end{equation}
		because $\el \ge C\ge \frac{9(k-1)kq}{\rho}$.
		The inequalities \eqref{eq:bound_I_3}, \eqref{eq:bound_I_4} and \eqref{eq:bound_I_5} show that we can apply Lemma \ref{lem:embedding_complete} to $G[N(A)\cap W^{(1)},\ldots,N(A)\cap W^{(k-1)}]$ with $r_{\ref{lem:embedding_complete}}=k-1$ and $q_{\ref{lem:embedding_complete}}=2q$.
	\end{proof}
	
	We also require another consequence of Lemma \ref{lem:embedding_complete}, stated below.
	
	\begin{lemma}\label{lem:bound_J}
		Given integers $k \ge 3$ and $q\ge 1$ and a real number $\eta \in \left(\frac{2q-1}{2(k-1)q},1\right)$, there exist an integer $C=C(k,q,\eta)$ and a positive $\de=\de(k,q,\eta)$ such that the following holds for every integer $\el \ge C$ and every $\eps \in (0,\de)$.
		Let $G=(V_1\cup\ldots\cup V_{\el},E)$ be a balanced $\el$-partite graph containing no copy of $K_{k-1}^{+q}(2q)$ whose vertices are in different parts of $G$. Assume $(X^{(s)}_i,Y^{(s)}_i)_{s\le k-1, i\le \el}$ are pairs of vertex sets satisfying:
		\begin{compactitem}
			\item[\rm (i)] For $i\le \el$ and $s\le k-1$, $Y_i^{(1)},\ldots, Y_i^{(k-1)}$ are disjoint subsets of $V_i$ and $X_i^{(s)}\subseteq Y_i^{(s)}$,
			\item[\rm (ii)] For $i\le \el$ and $s\le k-1$, $\card{X^{(s)}_i}= \left(\frac{1}{k-1}\pm\eps\right)\card{V_i}$,
			\item[\rm (iii)] For $s \le k-1$ and $v\in \bigcup_{i\le \el, t\ne s}X_i^{(t)}$, $\deg(v,\bigcup_{i\le \el}X_i^{(s)}) \ge \card{\bigcup_{i\le \el}X_i^{(s)}}-\eps \cdot v(G)$. 
		\end{compactitem}
		For $i\le \el$ and $s\le k-1$, let
		$B_i^{(s)}$ denote a subset of $Y_i^{(s)}$ consisting of all vertices $v$ with $\deg(v,\bigcup_{j\le \el}X_j^{(t)})<\eta\cdot v(G)$ for some $t \ne s$.
		For $s \le k-1$, write $\cM_{(s)}$ for a maximal matching in the induced subgraph $G\big[\bigcup_{i\le\el}Y_i^{(s)}\setminus B_i^{(s)}\big]$ of $G$ whose vertices are in different parts of $G$, and set $J=\{j\in [\el]: V_j \ \textrm{contains some vertex in $\bigcup_{s\le k-1} \cM_{(s)}$}\}$. Then, $\card{J}<2(k-1)q$.
	\end{lemma}
	\begin{proof}
		Choose
		\[
		C=\frac{4(k-2)}{\eta'} \ \text{and} \ \de=\min\biggl\{ \frac{q\eta'}{2q-1},\frac{\eta'}{4(k-2)}\biggl\}, \ \text{where $\eta'=\eta-\frac{2q-1}{2(k-1)q}$}.
		\]
		Notice that $\eta'>0$ as  $\eta \in \left(\frac{2q-1}{2(k-1)q},1\right)$. 
		We prove the statement by contradiction. Suppose that $\cM_{(s)}$ contains a matching $\{x_1x_2,\ldots,x_{2q-1}x_{2q}\}$ of size $q$ for some $s \le k-1$. 
		Let $X^{(t)}$ denote the vertex set $\bigcup_{i}X_i^{(s)}$ for $s\le k-1$. For $t \ne s$, define
		$W_{(t)}= \bigcup_{i}\bigl(N(x_1,\ldots,x_{2q}) \cap X_i^{(t)}\bigl)$. Then property (i) implies that $W^{(1)},\ldots,W^{(k-1)}$ are disjoint subsets of $V(G)$.
		We shall apply Lemma \ref{lem:embedding_complete} to find a copy of $K_{k-2}(2q)$ in $G \bigl[W_{(1)},\ldots,\widehat{W_{(s)}},\ldots,W_{(k-1)}\bigl]$ whose vertices are in different parts of $G$ (here $\widehat{W_{(s)}}$ stands for the empty set). Since this copy lies in $N(x_1,\ldots,x_{2q})$ and since $\{x_1x_2,\ldots,x_{2q-1}x_{2q}\}$ is a matching, $G$ contains a copy of $K_{k-1}^{+q}(2q)$ whose vertices belong to different parts of $G$, which is impossible. The remaining task is thus to verify the assumptions of Lemma \ref{lem:embedding_complete}. Indeed, from the definition of $W_{(t)}$ we see that, for $t\ne s$,
		\begin{equation}
		\label{eq:bound_J_1}
		W_{(t)}= \dot\bigcup_{i}\bigl(N(x_1,\ldots,x_{2q}) \cap X_i^{(t)}\bigl).
		\end{equation}
		By the definition of $\cM_{(s)}$, we have $\deg(x,X^{(t)}) \ge \eta \cdot v(G)$ for $x\in\{x_1,\ldots,x_{2q}\}$ and $t\ne s$. Hence
		\begin{align}\label{eq:bound_J_2}
		\notag	\card{W_{(t)}} &=\card{N(x_1,\ldots,x_{2q})\cap X^{(t)}}\ge 2q \eta\cdot v(G)-(2q-1)\card{X^{(t)}}\\
		& \hspace*{2cm}\overset{(ii)}{\ge} 2q\eta\cdot v(G)-(2q-1)\left(\tfrac{1}{k-1}+\eps\right)v(G) \ge q\eta'\cdot v(G)
		\end{align}
		for $\eps \le \de \le \frac{q\eta'}{2q-1}$.
		Together with the assumption $\el \ge C = \frac{4(k-2)}{\eta'}$, this inequality implies that, for $i\le \el$ and $t\ne s$,
		\begin{equation}
		\label{eq:bound_J_3}
		\card{N(x_1,\ldots,x_{2q})\cap X_i^{(t)}} \le \card{V_i}=\frac{v(G)}{\el} \le \tfrac{q\eta'}{4(k-2)q}\cdot v(G) \le \tfrac{1}{4(k-2)q}\cdot \card{W_{(t)}}.
		\end{equation}
		On the other hand, we can derive from property (iii) that, for $v\in \bigcup_{i\le \el,p\notin \{s,t\}}X_i^{(p)}$,
		\begin{equation}
		\label{eq:bound_J_4}
		\card{W_{(t)}}-\deg(v,W_{(t)})\le \eps \cdot v(G) \le \tfrac{q\eta'}{4(k-2)q} \cdot v(G) \overset{	\eqref{eq:bound_J_2}}{\le} \tfrac{1}{4(k-2)q}\cdot\card{W_{(t)}},
		\end{equation}
		assuming $\eps \le \de \le \frac{\eta'}{4(k-2)}$.
		It follows from \eqref{eq:bound_J_1}, \eqref{eq:bound_J_3} and \eqref{eq:bound_J_4} that we can apply Lemma \ref{lem:embedding_complete} to $G \bigl[W_{(1)},\ldots,\widehat{W_{(s)}},\ldots,W_{(k-1)}\bigl]$ with $r_{\ref{lem:embedding_complete}}=k-2$ and $q_{\ref{lem:embedding_complete}}=2q$. 
	\end{proof}
	
	We are now ready to prove Theorem \ref{thm:key}.
	
	\begin{proof}[Proof of Theorem \ref{thm:key}]
		Let $k=\chi(H)$. If $k=2$, then $H$ is a matching. The density condition implies that there is at least one edge between any two parts of $G$. Hence $G$ contains a  matching of size $ \frac{\el}{2}\ge e(H)$ whose vertices are in different parts of $G$.
		So from now on we can focus on the case when $k\ge 3$. Moreover, as discussed in Remark \ref{rmk:universal_graph}, we can suppose that $H=K^{+q}_{k-1}(2q)$ for some positive integer $q$. 
		To prove Theorem \ref{thm:key}, we assume to the contrary that $G$ does not contain a copy of $H$ whose vertices are in different parts of $G$.  
		Without loss of generality we can suppose that each part of $G$ has exactly $m$ vertices, where $m$ is a sufficiently large integer. 
		Otherwise, multiply each vertex in each part $V_i$ by a factor of $\frac{m}{\card{V_i}}$, which has no effect on the densities, and creates no copy of $H$ whose vertices lie in different parts of $G$.
		
		Choose $\el=\max\{C_{\ref{lem:induced_subgraph}}(k,q,\eps),1/\eps\}$, where $\eps>0$ is sufficiently small (to be specified later).
		Let $\el_1=\frac{\el}{2(k-1)!}$, 
		$\el_2=\el_1-(k-1)$, $\el_3=\frac{\el_2}{(k-1)!}$ and $\el_4=\el_3-2(k-1)q-D$, where $D=D_{\ref{lem:bound_I}}\left(k,q,\frac{1}{(6q+10)(k-1)(k-1)!}\right)$.
		Note that the parameters $\el$ and $\el_i$ both grow as $\Omega(1/\eps)$. 
		
		Our goal is to find an infracolourable struture in $G$. In the first step, we apply Lemma \ref{lem:induced_subgraph} to $G$ with $k_{\ref{lem:induced_subgraph}}=k$, $q_{\ref{lem:induced_subgraph}}=q$ and $\eps_{\ref{lem:induced_subgraph}}=\eps<\frac{1}{8k^2q}$ to obtain an induced $(k-1)$-colourable subgraph $F$ of $G$ whose vertex classes $X^{(1)}, \ldots, X^{(k-1)}$ satisfy 
		\begin{equation}
		\label{eq:key_1}
		\card{X^{(s)}}=
		\left(\tfrac{1}{k-1}\pm\eps\right)n \ \text{for $s \le k-1$},\\
		\end{equation}
		\begin{equation}
		\label{eq:key_2}
		\deg(v,X^{(s)}) \ge \card{X^{(s)}}-\eps n \ \text{for $s \le k-1$ and $v\in \bigcup_{t\ne s}X^{(t)}$}.
		\end{equation}
		
		Let $T=V(G) \setminus V(F)$. The inequality \eqref{eq:key_1} implies that $\card{T} \le k\eps n$. This forces $\card{T_i} \le 2k\eps m$ for at least half of indices $i \le \el$.
		Since $\el_1=\frac{\el}{2(k-1)!}$, by the pigeon hole principle we can relabel the $V_i$ and the $X^{(s)}$ such that
		$\card{X_i^{(1)}} \ge \card{X_i^{(2)}}\ge\ldots \ge \card{X_i^{(k-1)}}$ and
		\begin{equation}\label{eq:key_3}
		\card{T_i} \le 2k\eps m \ \text{for $i\le \el_1$}.
		\end{equation}
		Hence we can apply Lemma \ref{lem:balanced} with $\eps_{\ref{lem:balanced}}=2k\eps<\frac14$ to find a subset $I_0 \in \binom{\bN}{k-1}$ such that $\card{X_i^{(s)}}=\left(\tfrac{1}{k-1}\pm k\sqrt{2k\eps}\right)m$ for $s\le k-1$ and $i\in [\el_1]\setminus I_0$. By reordering parts if necessary, we may assume that
		\begin{equation}
		\label{eq:key_4}
		\card{X_i^{(s)}}=\left(\tfrac{1}{k-1}\pm k\sqrt{2k\eps}\right)m \quad \text{for $s\le k-1$ and $i\le \el_2$}.
		\end{equation} 
		
		For $i \le \el_2$ we shall partition $V_i$ into $k-1$ subsets $Y_i^{(1)},\ldots, Y_i^{(k-1)}$ as follows. A vertex $v \in V_i$ is assigned to $Y_i^{(s)}$ if
		$\deg\big(v,\bigcup_{j \le \el_2}X_j^{(s)}\big)= \min_{t\le k-1} \deg\big(v,\bigcup_{j \le \el_2}X_j^{(t)}\big)$; if there are more than one such index $s$, arbitrarily choose one of them.
		
		\begin{claim}\label{eq:key_5}
			$X_i^{(s)} \subseteq Y_i^{(s)} \subseteq X_i^{(s)}\dot\cup T_i$ and $\card{Y_i^{(s)}}=\left(\tfrac{1}{k-1}\pm 2k\sqrt{2k\eps}\right)m$ for $s \le k-1$ and $i \le \el_2$.
		\end{claim}
		\begin{proof}
			Let $v$ be an arbitrary vertex of $X_i^{(s)}$. Since $X^{(s)}$ is an independent set of $G$, $v$ has no neighbours in $\bigcup_{j \le \el_2}X_j^{(s)}$. It thus follows from the definition of $Y_i^{(s)}$ that $v \in Y_i^{(s)}$, and so $X_i^{(s)}$ is a subset of $Y_i^{(s)}$. Combining with the fact that $V_i=\left(\dot\bigcup_{s} X_i^{(s)}\right)\dot\cup T_i=\dot\bigcup_{s}Y_i^{(s)}$, we conclude that $Y_i^{(s)}\subseteq X_i^{(s)}\dot\cup T_i$ for $i \le \el_2$ and $s\le k-1$.
			
			As $X_i^{(s)}$ is a subset of $Y_i^{(s)}$, \eqref{eq:key_4} tells us that $\card{Y_i^{(s)}}\ge \card{X_i^{(s)}} \ge \left(\tfrac{1}{k-1}-k\sqrt{2k\eps}\right)m$ for $i \le \el_2$ and $s\le k-1$. Using \eqref{eq:key_3} and \eqref{eq:key_4}, we get 
			\[
			\card{Y_i^{(s)}}\le \card{X_i^{(s)}}+\card{T_i} \le \left(\tfrac{1}{k-1}+k\sqrt{2k\eps}+2k\eps\right)m\le \left(\tfrac{1}{k-1}+2k\sqrt{2k\eps}\right)m
			\]
			for $i\le \el_2$ and $s\le k-1$, where the first inequality holds since $Y_i^{(s)}$ is a subset of $X_i^{(s)}\cup T_i$. 
		\end{proof} 
		
		Let $I=\bigl\{i\in [\el_2]: \exists \ v_i \in V_i \ \text{with $\deg(v_i,X_1^{(s)}\cup\ldots\cup X_{\el_2}^{(s)}) \ge \tfrac{1}{(6q+10)(k-1)!} \cdot \frac{\el_2m}{k-1}$ for $s \le k-1$}\bigl\}$. We shall show that $I$ has bounded size.
		
		\begin{claim} \label{claim:bound_I} 
			$\card{I}\le D$.
		\end{claim}
		\begin{proof}
			We require $\eps$ to be small enough so that $\max\{k\sqrt{2k\eps},k^k\eps\} < \de_{\ref{lem:bound_I}}\left(k,q,\frac{1}{(6q+10)(k-1)(k-1)!}\right)$, and $\el_2 \ge C_{\ref{lem:bound_I}}\left(k,q,\frac{1}{(6q+10)(k-1)(k-1)!}\right)$.
			By \eqref{eq:key_4}, $\card{X_i^{(s)}}=\left(\tfrac{1}{k-1}\pm k\sqrt{2k\eps}\right)m$ for $s\le k-1$ and $i\le \el_2$. Moreover, for  $s\le k-1$ and $v\in \bigcup_{i\le \el_2, t\ne s}X_i^{(t)}$, we have 
			\begin{align*}
			\deg(v,\bigcup_{i\le \el_2}X_i^{(s)}) &\ge \card{\bigcup_{i\le \el_2}X_i^{(s)}}+\deg(v,X^{(s)})-\card{X^{(s)}}\\
			&\overset{\eqref{eq:key_2}}{\ge} \card{\bigcup_{i\le \el_2}X_i^{(s)}}-\eps n \ge \card{\bigcup_{i\le \el_2}X_i^{(s)}}-k^k\eps \el_2 m.
			\end{align*}
			Therefore, we can apply Lemma \ref{lem:bound_I} to $G[V_1\cup\ldots\cup V_{\el_2}]$ with input $k_{\ref{lem:bound_I}}=k$, $q_{\ref{lem:bound_I}}=q$ and $\eta_{\ref{lem:bound_I}}=\frac{1}{(6q+10)(k-1)(k-1)!}$ to conclude that $\card{I} \le D_{\ref{lem:bound_I}}\left(k,q,\frac{1}{(6q+10)(k-1)(k-1)!}\right)=D$. 
		\end{proof}
		
		As $\el_3=\frac{\el_2}{(k-1)!}$, by reordering the $V_i$ and $Y^{(s)}$ if necessary we can ensure
		\begin{equation}
		\label{eq:key_6}
		V_i=\dot\bigcup_{s}Y_i^{(s)} \ \text{and} \ \card{Y_i^{(1)}} \ge \card{Y_i^{(2)}}\ge \ldots \ge \card{Y_i^{(k-1)}} \quad \textrm{for $i \le \el_3$}.
		\end{equation}
		For $i \le \el_3$ and $s \le k-1$, let
		$B_i^{(s)}$ be the set of all vertices $v \in Y_i^{(s)}$ with the property that $\deg(v,X_1^{(t)}\cup\ldots\cup X_{\el_3}^{(t)})<\tfrac{2q}{2q+1}\cdot\frac{\el_3m}{k-1}$ for some $t \ne s$.
		For $s \le k-1$, let $\cM_{(s)}$ denote a maximal matching in $G\big[\bigcup_{i \le \el_3}Y_i^{(s)}\setminus B_i^{(s)}\big]$ whose vertices are in different parts of $G$, and write $J$ for the collection of all indices $j \in [\ell_3]$ so that $\bigcup_{s \le k-1} \cM_{(s)}$ contains some vertex in $V_j$.
		
		\begin{claim} \label{claim:bound_J} 
			$\card{J}<2(k-1)q$.
		\end{claim}
		\begin{proof}
			We shall apply Lemma \ref{lem:bound_J} to $G[V_1\cup\ldots V_{\el_3}]$ with $k_{\ref{lem:bound_J}}=k$, $q_{\ref{lem:bound_J}}=q$ and $\eta_{\ref{lem:bound_J}}=\frac{2q}{(k-1)(2q+1)}$ to get $\card{J}<2(k-1)q$. Note that $\card{X_i^{(s)}}=\left(\tfrac{1}{k-1}\pm k\sqrt{2k\eps}\right)m$ for $s\le k-1$ and $i\le \el_3$, by \eqref{eq:key_4}. Furthermore, for $s\le k-1$ and $v\in \bigcup_{i\le \el_3, t\ne s}X_i^{(t)}$, we have 
			\[
			\deg(v,\bigcup_{i\le \el_3}X_i^{(s)}) \overset{\eqref{eq:key_2}}{\ge} \card{\bigcup_{i\le \el_2}X_i^{(s)}}-\eps n \ge \card{\bigcup_{i\le \el_2}X_i^{(s)}}-k^{2k}\eps \el_3 m.
			\]
			Finally, we can choose $\eps$ sufficiently small so that $\max\{k\sqrt{2k\eps},k^{2k}\eps\}<\de_{\ref{lem:bound_J}}\left(k,q,\frac{2q}{(k-1)(2q+1)}\right)$ and $\el_3 \ge C_{\ref{lem:bound_J}}\left(k,q,\frac{2q}{(k-1)(2q+1)}\right)$.
		\end{proof}
		
		From Claims \ref{claim:bound_I} and \ref{claim:bound_J} we can assume (relabelling parts once more if necessary) that $\{1,\ldots,\el_3\}\setminus(I\cup J)=\{1,\ldots,\el_4\}$. For $i \le \el_4$ and $s \le k-1$, let 
		$D_i^{(s)}$ be the set consisting of all vertices $v \in Y_i^{(s)}$ such that $\deg(v,Y_1^{(t)}\cup\ldots\cup Y_{\el_4}^{(t)})<\tfrac{2q+1}{2q+2}\cdot \frac{\el_4m}{k-1}$ for some $t \ne s$.
		
		\begin{claim} \label{claim:find_infracolourable}
			The $\el_4$-partite graph
			$G[V_1\cup\ldots\cup V_{\el_4}]$ together with pairs $(D_i^{(s)},Y_i^{(s)})_{s\le k-1, i\le \el_4}$ of vertex sets form an $(\tfrac{1}{6q+9},k,\el_4)$-infracolourable structure.
		\end{claim}
		\begin{proof}
			We have to verify the following three properties:
			\begin{compactitem}
				\item[\rm (i)] For $i\le \el_4$, $V_i=\dot\bigcup_{s\le k-1} Y_i^{(s)}$ and $\card{Y_i^{(1)}}\ge \card{Y_i^{(2)}} \ge \ldots \ge \card{Y_i^{(k-1)}}$;
				\item[\rm (ii)] For $i\le \el_4$ and $s\le k-1$, $D_i^{(s)} \subseteq Y_i^{(s)}$ and $\bigcup_{i\le \el_4}Y_i^{(s)}\setminus D_i^{(s)}$ is an independent set;
				\item[\rm (iii)] For $s\le k-1$, every vertex $v\in \bigcup_{i\le \el_4}D_i^{(s)}$ has at most $\tfrac{1}{6q+9} \cdot \frac{\el_4m}{k-1}$ neighbours in $\bigcup_{i\le \el_4}Y_i^{(s)}$ and at least $\tfrac{1}{2q+3} \cdot \frac{\el_4m}{k-1}$ non-neighbours in $\bigcup_{i\le \el_4}V_i\setminus Y_i^{(s)}$.
			\end{compactitem}
			Property (i) follows directly from \eqref{eq:key_6}. For (ii), we observe that $B_i^{(s)} \subseteq D_i^{(s)}$ for $i \le \el_4$ and $s \le k-1$. We then deduce property (ii) from the maximality of $\cM_{(s)}$. 
			For (iii), we consider an arbitrary vertex $v \in \bigcup_{i\le \el_4}D_i^{(s)}$. Assume to the contrary that $\deg(v,\bigcup_{i \le \el_4}Y_i^{(s)})> \tfrac{1}{6q+9}\cdot\tfrac{\el_4m}{k-1}$. 
			Then, by Claim \ref{eq:key_5}, we obtain
			\[
			\deg(v,\bigcup_{i \le \el_4}X_i^{(s)}) \ge \deg(v,\bigcup_{i \le \el_4}Y_i^{(s)})-\card{\bigcup_{i\le \el_4}T_i} \\
			\overset{\eqref{eq:key_3}}{\ge} \tfrac{1}{6q+9}\cdot\tfrac{\el_4m}{k-1}-2k\eps \el_4 m > \tfrac{1}{(6q+10)(k-1)!}\cdot\tfrac{\el_2m}{k-1}
			\]
			for $\eps$ sufficiently small. On the other hand, by (ii), we must have $v\in \bigcup_{i\le \el_4}D_i^{(s)}\subseteq\bigcup_{i\le \el_4}Y_i^{(s)}$, and so $\deg(v,\bigcup_{i \le \el_2}X_i^{(t)}) \ge \deg(v,\bigcup_{i \le \el_2}X_i^{(s)})$ for all $t \le k-1$.
			Therefore,
			\begin{equation*}
			\deg(v,\bigcup_{i \le \el_2}X_i^{(t)}) \ge \deg(v,\bigcup_{i \le \el_2}X_i^{(s)})\ge \deg(v,\bigcup_{i \le \el_4}X_i^{(s)}) > \tfrac{1}{(6q+10)(k-1)!}\cdot\tfrac{\el_2m}{k-1}
			\]
			for $t\le k-1$, as $v\in \bigcup_{i\le \el_4}Y_i^{(s)}$. This contradicts the fact that $\{1,\ldots,\el_4\} \cap I=\emptyset$.
			Finally, by the definition of $\bigcup_{i\le \el_4}D_i^{(s)}$, there exists $t\ne s$ such that $\deg(v,\bigcup_{i \le \el_4}Y_i^{(t)})<\tfrac{2q+1}{2q+2}\cdot\frac{\el_4m}{k-1}$. Consequently, the number of non-neighbours of $v$ in in $\bigcup_{i \le \el_4}Y_i^{(t)}$ is at least
			\[
			\card{\bigcup_{i \le \el_4}Y_i^{(t)}}-\tfrac{2q+1}{2q+2}\cdot\tfrac{\el_4m}{k-1}\overset{\text{Claim \ref{eq:key_5}}}{\ge} \left(\tfrac{1}{k-1}-2k\sqrt{2k\eps}\right)\el_4m-\tfrac{2q+1}{2q+2}\cdot\tfrac{\el_4m}{k-1} >\tfrac{1}{2q+3}\cdot\tfrac{\el_4m}{k-1},
			\]
			assuming $\eps$ is sufficiently small.
			\end{proof}
			
			Claim \ref{claim:find_infracolourable} tells us that $G[V_1\cup\ldots\cup V_{\el_4}]$ is the base graph of an $(\tfrac{1}{6q+9},k,\el_4)$-infracolourable structure. By appealing to Lemma \ref{lem:infracolourable}, we can find two indices $1\le i<j\le \el_4$ with $d(V_i,V_j)\le \frac{k-2}{k-1}$, contradicting the assumption that $d(V_i,V_j)>\frac{k-2}{k-1}$. This completes our proof of Theorem \ref{thm:key}.
			\end{proof}

			\section{Proof of Theorem \ref{thm:complete_plus-weak}}
			\label{sect:complete_plus}
			In this section we shall prove a stronger version of Theorem \ref{thm:complete_plus-weak}.  
			
			\begin{theorem} \label{thm:complete_plus-strong}
			Let $k$ and $\el$ be integers with $k\ge 3$ and $\el \ge e^{2/c}$, where $c$ is a real number with $0<c \le k^{-(k+6)k}/2$. Suppose that $G=\left(V_1\cup \ldots \cup V_{\el},E\right)$ be a balanced $\el$-partite graph on $n$ vertices such that
			\[
			d(V_i,V_j) \ge \tfrac{k-2}{k-1} \quad \text{for $i \ne j$}.
			\] 
			Then, $G$ either contains a copy of $K_{k-1}^{+}\bigl(\lfloor c \ln n \rfloor, \ldots, \lfloor c \ln n \rfloor, \lfloor n^{1-2\sqrt{c}}\rfloor \bigl)$ or is isomorphic to a graph in $\mathcal{G}_{\el}^k$.
			\end{theorem}
			
			The idea of the proof is similar to that of Theorem \ref{thm:key}. We assume that $G$ does not contain a copy of $K_{k-1}^{+}\bigl(\lfloor c \ln n \rfloor, \ldots, \lfloor c \ln n \rfloor, \lfloor n^{1-2\sqrt{c}}\rfloor \bigl)$. We wish to show that $G$ is isomorphic to a graph in the family $\mathcal{G}_{\el}^k$. For this purpose, we apply the stability lemma (Lemma \ref{lem:induced-stability_complete}) to find an induced $(k-1)$-colourable subgraph of $G$ which almost spans $V(G)$. We then use the embedding lemma (Lemma \ref{lem:embedding_complete-plus}) showing that $G$ contains a large infracolourable structure. To conclude the proof, we shall use a bootstrapping argument (Lemma \ref{bootstrap}) which allows leveraging a weak structure result into a strong structure result. 
			
			In the proof of Theorem \ref{thm:complete_plus-strong} we shall need the following stability lemma.
			
			\begin{lemma}\label{lem:induced-stability_complete}
			Let $k$ and $\el$ be integers with $k\ge 3$ and $\el \ge e^{2/c}$, where $c$ is a real number with $0< c \le k^{-(k+6)k}/2$. Let $G=(V_1\cup\ldots\cup V_{\el},E)$ be a balanced $\el$-partite graph such that $d(V_i,V_j)\ge \frac{k-2}{k-1}$ for $i\ne j$. If $G$ does not contain a copy of $K_{k-1}^{+}\bigl(\lfloor c \ln n \rfloor, \ldots, \lfloor c \ln n \rfloor, \lfloor n^{1-2\sqrt{c}}\rfloor \bigl)$, then $G$ has an induced $(k-1)$-colourable subgraph $F$ whose vertex classes $X^{(1)}, \ldots, X^{(k-1)}$ satisfy the following properties with $\eps=4\el^{-1/2}$   
			\begin{compactitem}
			\item[\rm (i)] For $s \le k-1$, $\card{X^{(s)}}=\left(\tfrac{1}{k-1}\pm k\eps\right)n$;
			\item[\rm (ii)] For $s \le k-1$ and $v\in \bigcup_{t\ne s}X^{(t)}$, $\deg(v,X^{(s)})\ge \card{X^{(s)}}-k\eps n$. 
			\end{compactitem}
			\end{lemma}
			
			To prove the above statement we need a stability lemma of Nikiforov \cite[Theorem 3]{Nikiforov10}.
			
			\begin{lemma}  
			\label{lem:Nikiforov_stability}
			Let $k\ge 3$ be an integer, and let $c$ and $\de$ be positive real numbers with $c < k^{-(k+6)k}/2$ and $\de < \frac{1}{8k^8}$. Suppose that $G$ is a graph of order $n \ge e^{2/c}$ with $e(G) \ge \left(\frac{k-2}{k-1}-\de\right)\binom{n}{2}$. If $G$ has no copy of $K_{k-1}^{+}\bigl(\lfloor c \ln n \rfloor, \ldots, \lfloor c \ln n \rfloor, \lfloor n^{1-2\sqrt{c}}\rfloor \bigl)$, then
			$G$ contains an induced $(k-1)$-colourable subgraph $F$ of order $v(F)\ge(1-2\sqrt{\de})n$ and minimum degree $\de(F)\ge\left(\frac{k-2}{k-1}-4\sqrt{\de}\right)n$.  
			\end{lemma}
			
			\begin{proof}[Proof of Lemma \ref{lem:induced-stability_complete}]
			By the assumption, $\card{V_1}=\ldots=\card{V_{\el}}=\frac{n}{\el}:=m$. Together with the density condition, we conclude that $e(G) \ge \binom{\el}{2}\tfrac{k-2}{k-1}m^2\ge\left(\tfrac{k-2}{k-1}-\tfrac{1}{\el}\right)\tfrac{(\el m)^2}{2} = \left(\tfrac{k-2}{k-1}-\tfrac{1}{\el}\right)\tfrac{n^2}{2}$.
			Notice that $c \le k^{-(k+6)k}$, $\frac{1}{\el}<\frac{1}{8k^8}$ and $n\ge e^{2/c}$. 
			Thus, by applying Lemma \ref{lem:induced-stability_complete} to $G$ with $\de_{\ref{lem:induced-stability_complete}}=\frac{1}{\el}$ we obtain an $(k-1)$-colourable induced subgraph $F=G[X^{(1)}\cup\ldots\cup X^{(k-1)}]$ of $G$ with 
			$v(F)>(1-\eps)n$ and $\de(F)\ge\left(\tfrac{k-2}{k-1}-\eps\right)n$. Since $\de(F)\ge \left(\tfrac{k-2}{k-1}-\eps\right)n$ and since $X^{(s)}$ is an independent set, we must have 
			\[
			\card{X^{(s)}} \le n-\de(F)\le \left(\tfrac{1}{k-1}+\eps\right)n  
			\]
			for $s\le k-1$. This implies that
			\[
			\card{X^{(s)}} \ge v(F)-(k-2)\left(\tfrac{1}{k-1}+\eps\right)n \ge (1-\eps)n-(k-2)\left(\tfrac{1}{k-1}+\eps\right)n=\left(\tfrac{1}{k-1}-(k-1)\eps\right)n
			\]
			for $s\le k-1$. Therefore, for $s\le k-1$ and $v\in \bigcup_{t\ne s}X^{(t)}$, the number of non-neighbours of $v$ in $X^{(s)}$ is at most
			\[
			n-\card{X^{(s)}}-d_F(v) \le n-\left(\tfrac{1}{k-1}-(k-1)\eps\right)n-\left(\tfrac{k-2}{k-1}-\eps\right)n=k\eps n,
			\]
			as desired.
			\end{proof}
			
			The next ingredient we need is an embedding result.
			
			\begin{lemma} 
			\label{lem:embedding_complete-plus}
			Let $r \ge 2$ be an integer, and let $G$ be an $r$-colourable graph with vertex classes $W_{(1)},\ldots,W_{(r)}$ of the same size $h$. Suppose that $\deg(w,W_{(s)}) \ge \bigl(1-\frac{1}{r^2}\bigl)h$ for $s \le r$ and $w\in \bigcup_{t \ne s}W_{(t)}$. Then 
			\begin{compactitem}
			\item[\rm (1)] $G$ contains at least $\tfrac{1}{2}h^r$ copies of $K_r$,
			\item[\rm (2)] For every $\al \in (0,\tfrac{1}{4})$ and $s\le r$, $G$ contains a copy of  $K_r(\lfloor \al^r \ln h \rfloor,\ldots,\lfloor \al^r \ln h \rfloor, \lfloor h^{1-\al^{r-1}}\rfloor)$ whose $s$th vertex class is a subset of $W_{(s)}$.
			\end{compactitem}
			\end{lemma}
			
			The proof of the above lemma requires a simple result of Nikiforov \cite[Lemma 5]{Nikiforov10}. 
			
			\begin{lemma} \label{technical}
			Let $r \ge 2$ be an integer, and let $\al$ be a real number in $(0,\tfrac{1}{4})$. Suppose that $B[U,W]$ is a bipartite graph with $\card{U}=p$ and $\card{W}=q$. If $p \ge 4\lfloor \al^r \ln q \rfloor$ and $e(B[U,W]) \ge \tfrac{1}{2}pq$, then $B[U,W]$ contains the complete bipartite graph $K(a,b)$ with $a=\lfloor \al^r \ln q \rfloor$ and $b=\lfloor q^{1-\al^{r-1}}\rfloor$.
			\end{lemma}
			
			\begin{proof}[Proof of Lemma \ref{lem:embedding_complete-plus}]
			(1) Let $w_s \in W_{(s)}$ for $s=1,\ldots, r$. Observe that $\{w_1,\ldots,w_r\}$ forms a clique of $G$ if and only if $w_s \in N(w_1,\ldots,w_{s-1}) \cap W_{(s)}$ for $s=2,\ldots, r$. In addition, $\card{N(w_1,\ldots,w_{s-1}) \cap W_{(s)}} \ge h-(s-1)\cdot\frac{h}{r^2}$. Thus, we can bound the number of copies of $K_r$ in $G$ from below by 
			\[
			h^r \cdot\prod_{s=1}^{r}\bigl(1-\tfrac{s-1}{r^2}\bigl) \ge h^r\cdot\left(1-\sum_{s=1}^{r}\tfrac{s-1}{r^2}\right)=\tfrac{r+1}{2r}\cdot h^r >\tfrac{1}{2}h^r.
			\]
			
			(2) We proceed by induction on $r$. The base case $r=2$ follows from the first assertion and Lemma \ref{technical}. For the induction step, assume that $r>2$. The induction hypothesis implies that $G[W_{(1)}\cup \ldots \cup W_{(r-1)}]$ contains a copy of $K_{r-1}(m)$ with $m=\lfloor \al^{r-1} \ln h \rfloor$.
			Let $U$ denote a set of $m$ disjoint copies of $K_{r-1}$ in $K_{r-1}(m)$.
			Define a bipartite graph $B[U,W_{(r)}]$ with vertex classes $U$ and $W_{(r)}$, joining $R \in U$ to $w \in W_{(r)}$ if $R \cup \{w\}$ is a clique. We see that $\card{U}=m$ and $\card{W_{(r)}}=h$. Since $0<\al <1/4$, we have $m=\lfloor \al^{r-1} \ln h \rfloor \ge \lfloor 4 \al^r \ln h \rfloor \ge 4 \lfloor \al^r \ln h \rfloor$.
			Furthermore, every vertex of $U$ has at least $h-r\cdot\frac{h}{r^2} \ge h/2$ neighbours in $W_{(r)}$.
			Hence $e(B[U,W_{(r)}]) \ge mh/2$. The assertion then follows from the base case $r=2$.
			\end{proof}

			In order to find a large infracolourable structure in $G$ we shall use the following consequence of Lemma \ref{lem:embedding_complete-plus}.
			
			\begin{lemma}\label{lem:plus_bound_I}
			Let $k \ge 3$ and $\el \ge 2$ be integers, and let $\eps$ and $\al$ be positive real numbers with $\eps<10^{-2}k^{-k}$ and $\al<\tfrac{1}{4}$. Suppose that $G=(V_1\cup\ldots\cup V_{\el},E)$ is a balanced $\el$-partite graph containing no copy of $K_{k-1}^{+}\bigl(\lfloor\al^{k-1}\ln(p)\rfloor,\ldots,\lfloor\al^{k-1}\ln(p)\rfloor, \lfloor p^{1-\al^{k-2}} \rfloor\bigl)$, where $p=\tfrac{1}{16(k-1)(k-1)!}\cdot v(G)$. Assume that $(X^{(s)}_i)_{s\le k-1, i\le \el}$ are vertex sets so that
			\begin{compactitem}
			\item[\rm (i)] For $i\le \el$, $X_i^{(1)},\ldots, X_i^{(k-1)}$ are disjoint subsets of $V_i$;
			\item[\rm (ii)] For $s \le k-1$ and $v\in \bigcup_{i\le \el, t\ne s}X_i^{(t)}$, $\deg(v,\bigcup_{i\le \el}X_i^{(s)}) \ge \card{\bigcup_{i\le \el}X_i^{(s)}}-\eps \cdot v(G)$. 
			\end{compactitem}
			Then, there are no vertices $v \in V(G)$ such that $\deg(v,\bigcup_{i\le \el}X_i^{(s)}) \ge p$ for all $s \le k-1$.
			\end{lemma}
			\begin{proof}
			Suppose for the contradiction that there is $v \in V(G)$ with $\deg(v,\bigcup_{i}X_i^{(s)}) \ge p$ for all $s\le k-1$. Then,
			for $s \le k-1$ there exists a subset
			\begin{equation}
			\label{eq:plus_bound_I_1}
			W_{(s)} \subseteq N(v) \cap  \left(\bigcup_{i}X_i^{(s)}\right) \ \textrm{with $\card{W_{(s)}}=p$}.
			\end{equation}
			By property (i), $W_{(1)},\ldots,W_{(k-1)}$ are disjoint subsets of $V(G)$.
			On the other hand, property (ii) shows that for all $s \le k-1$ and $v \in \bigcup_{t\ne s}W_{(t)}$ one has
			\begin{align}
			\label{eq:plus_bound_I_2}	
			\notag \deg\bigl(v,W_{(s)}\bigl) &\ge \card{W_{(s)}}-\eps \cdot v(G) \\
			&\ge \card{W_{(s)}}-\frac{1}{(k-1)^2}\cdot \frac{1}{16(k-1)(k-1)!}\cdot v(G)=\left(1-\frac{1}{(k-1)^2}\right)\cdot\card{W_{(s)}},
			\end{align}
			as $\eps<10^{-2}k^{-k}$.
			Finally, it follows from \eqref{eq:plus_bound_I_1} and \eqref{eq:plus_bound_I_2} that we can apply Lemma \ref{lem:embedding_complete-plus}(2) to the graph $G[W_{(1)}, \ldots,W_{(k-1)}]$ with $r_{\ref{lem:embedding_complete-plus}}=k-1$, $h_{\ref{lem:embedding_complete-plus}}=p$ and $\al_{\ref{lem:embedding_complete-plus}}=\al$ to find a copy of $K_{k-1}\bigl(\lfloor\al^{k-1}\ln(p)\rfloor,\ldots,\lfloor\al^{k-1}\ln(p)\rfloor, \lfloor p^{1-\al^{k-2}} \rfloor\bigl)$. Since $W_{(1)} \cup \ldots \cup W_{(k-1)}$ lies in the neighbour of $v$, $G$ contains a copy of $K_{k-1}^{+}\bigl(\lfloor\al^{k-1}\ln(p)\rfloor,\ldots,\lfloor\al^{k-1}\ln(p)\rfloor, \lfloor p^{1-\al^{k-2}} \rfloor\bigl)$, which contradicts our assumption. 
			\end{proof}
			
			To find a large infracolourable structure in $G$ we also require the following consequence of Lemma \ref{lem:embedding_complete-plus}.
			
			\begin{lemma}\label{lem:plus_bound_J}
			Let $k\ge 3$ and $\el \ge 2$ be integers, and let $\eps$ and $\al$ be positive real numbers with $\eps<\frac{1}{12k^3}$ and $\al<\tfrac{1}{4}$. Let $G=(V_1\cup\ldots\cup V_{\el},E)$ be a balanced $\el$-partite graph containing no copy of $K_{k-1}^{+}\bigl(\lfloor\al^{k-1}\ln(p)\rfloor,\ldots,\lfloor\al^{k-1}\ln(p)\rfloor, \lfloor p^{1-\al^{k-2}} \rfloor\bigl)$, where $p=\tfrac{1}{4(k-1)}\cdot v(G)$.
			Suppose $(X^{(s)}_i,Y^{(s)}_i)_{s\le k-1, i\le \el}$ are pairs of vertex sets which satisfy
			\begin{compactitem}
			\item[\rm (i)] For every $i \le \el$ and $s\le k-1$, $Y_i^{(1)},\ldots, Y_i^{(k-1)}$ are disjoint subsets of $V_i$ and $X_i^{(s)}\subseteq Y_i^{(s)}$;
			\item[\rm (ii)] For $i\le \el$ and $s\le k-1$, $ \card{X_i^{(s)}} = \left(\frac{1}{k-1}\pm \eps\right)\card{V_i}$;
			\item[\rm (iii)] For $s \le k-1$ and $v\in \bigcup_{i\le \el, t\ne s}X_i^{(t)}$, $\deg(v,\bigcup_{i\le \el}X_i^{(s)}) \ge \card{\bigcup_{i\le \el}X_i^{(s)}}-\eps \cdot v(G)$. 
			\end{compactitem}
			For $i\le \el$ and $s\le k-1$, let $B_i^{(s)}$ stands for a subset of $Y_i^{(s)}$ consisting of all vertices $v$ with $\deg(v,\bigcup_{j\le \el}X_j^{(t)})< \tfrac{2}{3(k-1)}\cdot v(G)$ for some $t \ne s$.
			Then, for $s\le k-1$, $\bigcup_{i\le \el}Y_i^{(s)}\setminus B_i^{(s)}$ is an independent set of $G$.
			\end{lemma}
			\begin{proof}
			We prove by contradiction. Suppose that there exists an edge $\{x,y\}\in E(G)$ with $x, y \in \bigcup_{i}Y_i^{(s)} \setminus B_i^{(s)}$. 
			Let $t\ne s$. By the definition of $\bigcup_{i}B_i^{(s)}$, both $\deg(x, \bigcup_{i}X_i^{(t)})$ and $\deg(y, \bigcup_{i}X_i^{(t)})$ are at least $\tfrac{2}{3(k-1)}\cdot v(G)$.
			Hence
			\begin{align*}
			\card{N(x,y) \cap \bigcup_{i}X_i^{(t)}} &\ge \deg(x, \bigcup_{i}X_i^{(t)})+\deg(y, \bigcup_{i}X_i^{(t)})-\card{\bigcup_{i}X_i^{(t)}}\\
			& \overset{(ii)}{\ge} \frac{4}{3(k-1)}\cdot v(G)-\left(\frac{1}{k-1}+\eps\right)\cdot v(G) \ge \frac{1}{4(k-1)}\cdot v(G),
			\end{align*}
			as $\eps<\frac{1}{12k^3}$. It means that there is a subset
			\begin{equation*}
			\label{eq:plus_10}
			W_{(t)} \subseteq N(x,y) \cap \bigcup_{i \le \el_3}X_i^{(t)} \ \textrm{with $\card{W_{(t)}}=\tfrac{1}{4(k-1)}\cdot v(G)$}.
			\end{equation*}
			On the other hand, it follows from property (ii) that
			$\card{\bigcup_{i}X_i^{(s)}}\ge \left(\tfrac{1}{k-1}-\eps\right)v(G) >\tfrac{1}{4(k-1)}\cdot v(G)$ for $0<\eps<\frac{1}{12k^3}$, and so there exists a subset 
			\begin{equation*}
			\label{eq:plus_11}
			W_{(s)} \subseteq \bigcup_{i}X_i^{(s)} \ \textrm{with $\card{W_{(s)}} =\tfrac{1}{4(k-1)}\cdot v(G)$}.
			\end{equation*}
			Analysis similar to that in the proof of Lemma \ref{lem:plus_bound_I} shows that $G[W_{(1)}, \ldots, W_{(k-1)}]$ must contain a copy of $K_{k-1}\bigl(\lfloor\al^{k-1}\ln(p)\rfloor,\ldots,\lfloor\al^{k-1}\ln(p)\rfloor, \lfloor p^{1-\al^{k-2}} \rfloor\bigl)$ whose $s$th vertex class is of size $\lfloor\al^{k-1}\ln(p)\rfloor$. Adding back vertices $x$ and $y$ to this class one gets a supgraph of the graph $K_{k-1}^{+}\bigl(\lfloor\al^{k-1}\ln(p)\rfloor,\ldots,\lfloor\al^{k-1}\ln(p)\rfloor, \lfloor p^{1-\al^{k-2}} \rfloor\bigl)$, contradicting the hypothesis. 
		\end{proof}
		
		The last component of the proof is a bootstrapping argument which allows us to leverage a weak structure result into a strong structure result. Roughly speaking, it says that if $G$ contains an $\tilde{\el}$-partite subgraph which is in $\cG_{\tilde{\el}}^k$, then $G$ must belong to $\cG_{\el}^k$.
		
		\begin{lemma}\label{bootstrap}
			Let $k\ge 3$ be an integer, and let $G=(V_1\cup\ldots\cup V_{\el},E)$ be an $\el$-partite graph with $\card{V_1}=\ldots=\card{V_{\el}}=m$ and $d(V_i,V_j)\ge \frac{k-2}{k-1}$ for all $i\ne j$. Suppose that there exist an integer $\tilde{\el}$ and disjoint subsets $Y_i^{(1)},\ldots,Y_i^{(k-1)}$ of $V_i$ for $1\le i\le \tilde{\el}$ so that $\card{Y_i^{(s)}}=\tfrac{m}{k-1}$ and $d(Y_i^{(s)},Y_j^{(t)})=1$ for all $i\ne j$ and $s\ne t$. If $G$ does not contain a copy of $K_{k-1}^{+}\big(\frac{\tilde{\el}m}{32k^2}\big)$, then $G$ is isomorphic to a graph in the family $\cG_{\el}^k$. 
		\end{lemma}
		
		\begin{proof}
		We wish to show that $G$ is isomorphic to a graph in $\cG_{\el}^k$. 
			According to Lemma \ref{lem:extremal_graphs}, it suffices to prove $G$ is $(k-1)$-colourable. 
			By the assumption, we have
			\begin{equation}
			\label{eq:bootstrap_1}
			\card{Y_i^{(s)}}=\tfrac{m}{k-1}, \ \textrm{$d(Y_i^{(s)},Y_j^{(t)})=1$ for $s\ne t$ and $1 \le i < j \le \tilde{\el}$}. 
			\end{equation}
			
			We shall show that for $v \in V(G) \setminus \left(V_1\cup \ldots \cup V_{\tilde{\el}}\right)$ there does not exist $s \le k-1$ with 
			\begin{equation}
			\label{eq:bootstrap_2}
			\deg\bigl(v,Y_1^{(s)} \cup \ldots \cup Y_{\tilde{\el}}^{(s)}\bigl)\ge 1, \ \textrm{$\deg\bigl(v,Y_1^{(t)} \cup \ldots \cup Y_{\tilde{\el}}^{(t)}\bigl) \ge \frac{\tilde{\el}m}{2k}$ for all $t \ne s$}. 
			\end{equation}
			We prove by contradiction. Suppose that \eqref{eq:bootstrap_2} holds. We can pick an index $i_0\in \{1,2,\ldots,\tilde{\el}\}$ with $N(v) \cap Y_{i_0}^{(s)} \ne \emptyset$ whose existence is guaranteed by \eqref{eq:bootstrap_2}. We then arbitrarily add other indices to get a subset $I_{(s)} \subset \{1,\ldots,\tilde{\el}\}$ of size $\frac{\tilde{\el}}{8k}$. It follows from \eqref{eq:bootstrap_1} and \eqref{eq:bootstrap_2} that for each $t \ne s$, there are at least $\frac{\tilde{\el}}{4}$ indices $i \le \tilde{\el}$ with $\deg(v,Y_i^{(t)}) \ge \frac{m}{4k}$. 
			Hence we can find $k-1$ disjoint subsets $I_{(1)}, \ldots, I_{(k-1)}$ of size $\frac{\tilde{\el}}{8k}$ of $\{1,\ldots,\tilde{\el}\}$ with the property that $\deg(v,Y_i^{(t)}) \ge \frac{m}{4k}$ for all $t \ne s$ and $i \in I_{(t)}$. By \eqref{eq:bootstrap_1}, $G\big[\bigcup_{i \in I_{(1)}}Y_i^{(1)},\ldots,\bigcup_{i \in I_{(k-1)}}Y_i^{(k-1)}\big]$ is a complete $(k-1)$-partite graph. In addition, we have $\card{N(v)\cap\bigcup_{i \in I_{(s)}}Y_i^{(s)}} \ge \card{N(v)\cap Y_{i_0}^{(s)}}>0$ and
			\begin{equation*}
			\card{N(v)\cap\bigcup_{i \in I_{(t)}}Y_i^{(t)}}=\sum_{i\in I_{(t)}}\deg(v,Y_i^{(t)})\ge \card{I_{(t)}}\cdot \frac{m}{4k}=\frac{\tilde{\el}m}{32k^2} \quad \text{for $t\ne s$}.
			\end{equation*}
			Therefore, by adding $v$ to the $s$th part of $G\big[\bigcup_{i \in I_{(1)}}Y_i^{(1)},\ldots,\bigcup_{i \in I_{(k-1)}}Y_i^{(k-1)}\big]$ we get a supergraph of $K_{k-1}^{+}\big(\frac{\tilde{\el}m}{32k^2}\big)$ in $G$, contradicting our assumption.
			
			We can infer from \eqref{eq:bootstrap_2} that 
			$\deg(v,V_1 \cup \ldots \cup V_{\tilde{\el}}) \le \tfrac{k-2}{k-1} \cdot \tilde{\el}m$ for all $v \in V(G) \setminus \left(V_1\cup \ldots \cup V_{\tilde{\el}}\right)$. By the density condition, equality must hold. Again \eqref{eq:bootstrap_2} shows that for each $v \in V(G) \setminus \left(V_1\cup \ldots \cup V_{\tilde{\el}}\right)$, 
			\begin{equation}\label{eq:bootstrap_3}
	N(v) \cap \bigl(V_1\cup\ldots\cup V_{\tilde{\el}}\bigl) = \bigcup_{i \le \tilde{\el}} V_i \setminus Y_i^{(s)} \quad \text{for some $s\le k-1$.}		
			\end{equation} If $v \in V_i$ for some $i>\tilde{\el}$, then we assign $v$ to $Z_i^{(s)}$. For $i\le \tilde{\el}$ we let $Z_i^{(s)}=Y_i^{(s)}$ for $s\le k-1$. If we denote $Z^{(s)}= \dot\bigcup_{i} Z_i^{(s)}$ for $s\le k-1$, then $V=\dot\bigcup_{s} Z^{(s)}$. To prove $G$ is $(k-1)$-colourable, it is enough to show that $Z^{(1)},\ldots,Z^{(k-1)}$ are independent sets. Suppose to the contrary that for some $s\le k-1$, $Z^{(s)}$ contains an edge $\{u,v\}$ with $u\in Z_{i_1}^{(s)}$ and $v\in Z_{i_2}^{(s)}$. We can easily find $k-1$ disjoint subsets $J_{(1)},\ldots,J_{(k-1)}$ of size $\frac{\tilde{\el}}{2(k-1)}$ of $[\tilde{\el}]\setminus \{i_1,i_2\}$. Let $W^{(s)}=\{u,v\}\cup \left(\bigcup_{i\in J_{(s)}}Y_i^{(s)}\right)$ and $W^{(t)}=\bigcup_{i\in J_{(t)}}Y_i^{(t)}$ for $t\ne s$. It follows from \eqref{eq:bootstrap_1} and \eqref{eq:bootstrap_3} that $G[W^{(1)},\ldots,W^{(k-1)}]$ is a complete $(k-1)$-colourable graph with $\card{W^{(t)}}\ge \frac{\tilde{\el}}{2(k-1)}\cdot \frac{m}{k-1}>\frac{\tilde{\el}m}{32k^2}$ for $t \le k-1$. Combining this with the assumption that $\{u,v\}\in E(G)$, we conclude that $G$ contains a copy of $K_{k-1}^{+}\big(\frac{\tilde{\el}m}{32k^2}\big)$, a contradiction.
		\end{proof}
		
		We now have all the necessary tools to prove Theorem \ref{thm:complete_plus-strong}.
		
		\begin{proof}[Proof of Theorem \ref{thm:complete_plus-strong}]
			For convenience, we write $H=K_{k-1}^{+}\big(\lfloor c \ln n \rfloor,\ldots,\lfloor c \ln n \rfloor, \lfloor n^{1-2\sqrt{c}}\rfloor \big)$ and
			$H^{-}=K_{k-1}\big(\lfloor c \ln n \rfloor,\ldots,\lfloor c \ln n \rfloor, \lfloor n^{1-2\sqrt{c}}\rfloor \big)$. Suppose $G$ has no copy of $H$. We wish to show that $G$ is isomorphic to a graph in $\cG_{\el}^k$. Since $G$ is a balanced $\el$-partite graph on $n$ vertices, each partition set of $G$ has size $n/\el:=m$. Let $\eps=4\el^{-1/2}$, $\el_1=\frac{\el}{2(k-1)!}-(k-1)$, $\el_2=\frac{\el_2}{(k-1)!}$ and $\el_3=\el_2-1$.
			
			By Lemma \ref{lem:induced-stability_complete}, $G$ must contain an induced $(k-1)$-colourable subgraph $F$ whose vertex classes $X^{(1)}, \ldots, X^{(k-1)}$ satisfy 
			\begin{equation}
			\label{eq:plus_1}
			\card{X^{(s)}}=\left(\tfrac{1}{k-1}\pm k\eps\right)n \ \text{for $s \le k-1$},\\
			\end{equation}
			\begin{equation}
			\label{eq:plus_2}
			\deg(v,X^{(s)}) \ge \card{X^{(s)}}-k\eps n \ \text{for $s \le k-1$ and $v\in \bigcup_{t\ne s}X^{(t)}$}.
			\end{equation}
			
   		    Let $T=V(G) \setminus V(F)$. As in the proof of Theorem \ref{thm:key}, by relabelling parts we can assume that
			\begin{equation}
			\label{eq:plus_4}
			\card{T_i} \le 2k^2\eps m, \ \text{and} \ \card{X_i^{(s)}}=\left(\tfrac{1}{k-1}\pm 2k^2\sqrt{\eps}\right)m \quad \text{for $i\le \el_1$ and $s\le k-1$}.
			\end{equation} 
			For $i \le \el_1$ we shall partition $V_i$ into $k-1$ subsets as follows. A vertex $v \in V_i$ is assigned to $Y_i^{(s)}$ if
			$\deg\big(v,\bigcup_{j \le \el_1}X_j^{(s)}\big)= \min_{t\le k-1}\deg\big(v,\bigcup_{j \le \el_1}X_j^{(t)}\big)$; if there are more than one such index $s$, arbitrarily pick one of them. 
			
			\begin{claim}\label{eq:plus_5}
				$X_i^{(s)} \subseteq Y_i^{(s)} \subseteq X_i^{(s)}\dot\cup T_i$ for  $i\le \el_1$ and $s\le k-1$.	
			\end{claim}
			\begin{proof}
				Because $X^{(s)}$ is an independent set in $G$, every vertex in $X_i^{(s)}$ has no neighbours in $\bigcup_{j \le \el_1}X_j^{(s)}$, and so $X_i^{(s)}$ is a subset of $Y_i^{(s)}$. Since $V_i=\left(\dot\bigcup_{s} X_i^{(s)}\right)\dot\cup T_i=\dot\bigcup_{s}Y_i^{(s)}$ and $X_i^{(s)}\subseteq Y_i^{(s)}$ for $i\le \el_1$ and $s \le k-1$, the inclusion relation $Y_i^{(s)} \subseteq X_i^{(s)}\dot\cup T_i$ holds for $i\le \el_1$ and $s\le k-1$. 	
			\end{proof}
			
			We proceed by showing that $\bigcup_{i\le \el_1}V_i$ does not contain a vertex which has relatively large degree to $\bigcup_{i \le \el_1}Y_i^{(s)}$ for all $s \le k-1$.  
			\begin{claim} \label{claim:plus_bound_I} There are no vertices $v \in \bigcup_{i\le \el_1}V_i$ with $\deg(v,\bigcup_{i \le \el_1}Y_i^{(s)}) \ge \tfrac{1}{15(k-1)(k-1)!}\cdot \el_1 m$ for all $s \le k-1$.
			\end{claim}
			\begin{proof}
				We can derive from \eqref{eq:plus_2} that, for $s\le k-1$ and $v\in \bigcup_{i\le \el_1,t \ne s}X_i^{(t)}$,
				\begin{align*}
				\deg(v,\bigcup_{i \le \el_1}X_i^{(s)})&\ge \card{\bigcup_{i \le \el_1}X_i^{(s)}}+\deg(v,X^{(s)})-\card{X^{(s)}}\\
				&\ge \card{\bigcup_{i \le \el_1}X_i^{(s)}}-k\eps n \ge \card{\bigcup_{i \le \el_1}X_i^{(s)}}-k^k\eps \cdot \el_1m.
				\end{align*}
				Applying Lemma \ref{lem:plus_bound_I} to $G[V_1\cup\ldots\cup V_{\el_1}]$ with  $k_{\ref{lem:plus_bound_I}}=k$, $\eps_{\ref{lem:plus_bound_I}}=k^k\eps$ and $\al_{\ref{lem:plus_bound_I}}=(2c)^{1/(k-1)}$, we conclude that either $G[V_1\cup\ldots\cup V_{\el_1}]$ contains a copy of $K_{k-1}^{+}\bigl(\lfloor\al^{k-1}\ln(p)\rfloor,\ldots,\lfloor\al^{k-1}\ln(p)\rfloor, \lfloor p^{1-\al^{k-2}} \rfloor\bigl)$ or there are no vertices $v\in V_1\cup\ldots\cup V_{\el_1}$ with $\deg(v,X_1^{(s)}\cup\ldots\cup X_{\el}^{(s)}) \ge p$ for all $s$, where $p=\tfrac{1}{16(k-1)(k-1)!}\cdot \el_1m$. Since $\al^{k-1}\ln(p)>c\ln(n)$, $p^{1-\al^{k-2}}>n^{1-2\sqrt{c}}$ and since $G$ has no copy of $K_{k-1}^{+}\bigl(\lfloor c \ln n \rfloor, \ldots, \lfloor c \ln n \rfloor, \lfloor n^{1-2\sqrt{c}}\rfloor \bigl)$, the former case is ruled out. The later case implies our statement. 
			\end{proof}
			
			Since $\el_2=\frac{\el_1}{(k-1)!}$, by reordering parts if necessary we can assume that 
			\begin{equation}
			\label{eq:plus_6}
			\card{Y_i^{(1)}} \ge \card{Y_i^{(2)}} \ge \ldots \ge \card{Y_i^{(k-1)}} \ \textrm{for $i \le \el_2$}. 
			\end{equation}
			For $i \le \el_2$ and $s \le k-1$, let us denote
			\[
			D_i^{(s)}=\biggl\{v \in Y_i^{(s)}: \deg(v,Y_1^{(t)}\cup\ldots\cup Y_{\el_2}^{(t)})< \tfrac{3}{4(k-1)}\cdot\el_2m \ \text{for some $t \ne s$}\biggl\}.
			\]
			
			\begin{claim} \label{claim:plus_bound_J}
				The vertex set $\bigcup_{i\le \el_2}Y_i^{(s)}\setminus D_i^{(s)}$ is an independent set of $G$ for $s \le k-1$.
			\end{claim}
			\begin{proof}
				For $i\le \el_2$ and $s\le k-1$, let $B_i^{(s)}$ be the vertex set consisting of all vertices $v \in Y_i^{(s)}$ such that $\deg(v,\bigcup_{i\le \el_2}X_i^{(t)})< \tfrac{2}{3(k-1)}\cdot \el_2m$ for some $t \ne s$.
				Note that, for $s\le k-1$ and $v\in \bigcup_{i\le \el_2,t\ne s}X_i^{(s)}$, one has
				\begin{align*}
				\deg(v,\bigcup_{i \le \el_2}X_i^{(s)})&\ge \card{\bigcup_{i \le \el_2}X_i^{(s)}}+\deg(v,X^{(s)})-\card{X^{(s)}}\\
				&\overset{\eqref{eq:plus_2}}{\ge} \card{\bigcup_{i \le \el_2}X_i^{(s)}}-k\eps n \ge \card{\bigcup_{i \le \el_2}X_i^{(s)}}-k^{2k}\eps \cdot \el_2m.
				\end{align*}
				This estimate together with Claim \ref{eq:plus_5} and \eqref{eq:plus_4} show  that we can apply Lemma \ref{lem:plus_bound_J} to $G[V_1\cup\ldots\cup V_{\el_2}]$ with $k_{\ref{lem:plus_bound_J}}=k$, $\eps_{\ref{lem:plus_bound_J}}=\max\{2k^2\sqrt{\eps},k^{2k}\eps\}$ and $\al_{\ref{lem:plus_bound_J}}=(2c)^{1/(k-1)}:=\al$ to conclude that either $G[V_1\cup\ldots\cup V_{\el_2}]$ contains  $K_{k-1}^{+}\bigl(\lfloor\al^{k-1}\ln(p)\rfloor,\ldots,\lfloor\al^{k-1}\ln(p)\rfloor, \lfloor p^{1-\al^{k-2}} \rfloor\bigl)$ or $\bigcup_{i\le \el_2}Y_i^{(s)}\setminus B_i^{(s)}$ is an independent set of $G$ for $s \le k-1$, where $p=\tfrac{1}{4(k-1)}\cdot\el_2m$. Since $\al^{k-1}\ln(p)>c\ln(n)$, $p^{1-\al^{k-2}}>n^{1-2\sqrt{c}}$ and since $G$ has no copy of $K_{k-1}^{+}\bigl(\lfloor c \ln n \rfloor, \ldots, \lfloor c \ln n \rfloor, \lfloor n^{1-2\sqrt{c}}\rfloor \bigl)$, the former case is ruled out. We can see that the later case implies our statement.
			\end{proof}
			
			Now we can find a large infracolourable structure in $G$, and then use Lemma \ref{bootstrap} to show that $G$ is isomorphic to a graph in $\cG_{\el}^k$.
			
			\begin{claim} 
				$G$ is isomorphic to a graph in the family $\cG_{\el}^k$.
			\end{claim}
			\begin{proof}
				Analogously to the proof of Claim \ref{claim:find_infracolourable}, we can infer from Claims \ref{claim:plus_bound_I} and \ref{claim:plus_bound_J}, \eqref{eq:plus_4} and \eqref{eq:plus_6} that ${G[V_1\cup\ldots\cup V_{\el_2}]}$ together with pairs $(D_i^{(s)},Y_i^{(s)})_{s\le k-1, i\le \el_2}$ form a $(\tfrac{1}{15},k,\el_2)$-infracolourable structure. By Lemma \ref{lem:infracolourable} this implies that $e(G[V_1\cup\ldots\cup V_{\el_2}]) \le \binom{\el_2}{2}\frac{k-2}{k-1}m^2$ and hence the equality must occur by the density condition. Appealing to Lemma \ref{lem:infracolourable} once again, we see that there exists $i_0\in \{0,1,\ldots,\el_2\}$ with $\card{Y_i^{(s)}}=\tfrac{m}{k-1}$ for all $s$ and all $i\in [\el_2]\setminus\{i_0\}$, and $d(Y_i^{(s)},Y_j^{(t)})=1$ for all $s\ne t$ and $1\le i<j \le \el_2$. Hence we can apply Lemma \ref{bootstrap} with $\tilde{\el}=\el_2-1$ to conclude that either $G$ contains a copy of $K_{k-1}^{+}\big(\frac{(\el_2-1)m}{32k^2}\big)$ or $G$ is isomorphic to a graph in $\cG_{\el}^k$. The former can not happen since $G$ has no copy of $K_{k-1}^{+}\big(\lfloor c \ln n \rfloor,\ldots,\lfloor c \ln n \rfloor, \lfloor n^{1-2\sqrt{c}}\rfloor \big)$ and since $\frac{(\el_2-1)m}{32k^2}>\max\{n^{1-2\sqrt{c}},c\ln n\}$. So $G$ must isomorphic to a graph in the family $\cG_{\el}^k$.
			\end{proof}
			
			This concludes our proof of Theorem \ref{thm:key}.
		\end{proof}

		\section{Missing proofs}\label{sect:missing_proofs}
		
		\subsection{Proof of Theorem \ref{thm:joints}} 
		\label{sect:joints}
		In this section we sketch a proof of Theorem \ref{thm:joints}.
		We follow essentially the proof of Theorem \ref{thm:complete_plus-strong}. We make the following alterations. Instead of Lemma \ref{lem:Nikiforov_stability} we use a stability result due to Bollob\'as and Nikiforov \cite[Theorem 9]{BN08}.
		\begin{lemma}
			\label{lem:stability_joints}
			Let $k \ge 2$ be an integer, and let $\de$ be a positive with $\de < \frac{1}{16k^8}$. Suppose that $G$ is a graph with $n > k^8$ vertices and $e(G) \ge \left(\frac{k-2}{k-1}-\de \right)\binom{n}{2}$ edges. Then, either
			$G$ contains a family of $k^{-(k+5)}n^{k-2}$ copies of $K_k$ sharing a common edge, or $G$ contains an induced $(k-1)$-colourable subgraph $F$ of size $v(F)\ge (1-2\sqrt{\de})n$ and minimum degree $\de(F)\ge \left(\frac{k-2}{k-1}-4\sqrt{\de}\right)n$.  
		\end{lemma}
		
		We replace Lemma \ref{lem:embedding_complete-plus} by the following embedding result.
		
		\begin{lemma}
			\label{lem:embedding_joints}
			Let $r \ge 2$ be an integer, and let $G$ be an $r$-colourable graph with classes $W_{(1)},\ldots,W_{(r)}$ of the same size $h$. Suppose that $\deg(v,W_{(s)}) \ge \bigl(1-\frac{1}{r^2}\bigl)h$ for $s \le r$ and $v \in \bigcup_{t \ne s}W_{(t)}$. Then for every pair $(s,t)$ with $s\ne t$, there is an edge between $W_{(s)}$ and $W_{(t)}$ which is contained in $\tfrac{1}{2}h^{r-2}$ copies of $K_r$. 
		\end{lemma}
		
		\begin{proof}
			According to Lemma \ref{lem:embedding_complete-plus}, $G$ contains at least $\tfrac{1}{2}h^r$ copies of $K_r$. Hence there exists an edge between $W_{(s)}$ and $W_{(t)}$ which is shared by at least $h^r/(2h^2)=\tfrac{1}{2}h^{r-2}$ copies of $K_r$.
		\end{proof}
		
		The remainder of the proof is similar to that of Theorem \ref{thm:complete_plus-strong}.
		
		\subsection{Proofs of Proposition \ref{prop:Erdos-Stone_stability} and Lemma \ref{lem:common_neighbours}}
		To prove Proposition \ref{prop:Erdos-Stone_stability} we shall require the Erd\H{o}s-Simonovits stability theorem (Erd\H{o}s \cite{Erdos67} and Simonovits \cite[Theorem 8]{Simonovits66}, and the graph removal lemma (Ruzsa and Szemer\'edi \cite{RS78}).
		
		\begin{theorem} [\bf{Stability theorem}]
			\label{stability theorem}
			For every graph $H$ and every $\eps>0$, there exist positive constants $\de=\de(H,\eps)$ and $C=C(H,\eps)$ so that the following holds for every integer $n\ge C$. Every $n$-vertex $H$-free graph with at least $\left(\frac{\chi(H)-2}{\chi(H)-1}-\de\right)\binom{n}{2}$ edges contains a $(\chi(H)-1)$-colourable subgraph of order at least $(1-\eps)n$ and minimum degree at least $\left(\frac{\chi(H)-2}{\chi(H)-1}-\eps\right)n$.
		\end{theorem}
		
		\begin{theorem} [{\bf Graph removal lemma}]
			\label{graph removal lemma}
			For every graph $H$ and every $\de>0$, there exists a positive constant $\ga=\ga(H,\de)$ such that every graph on $n$ vertices with at most $\ga n^{v(H)}$ copies of $H$ can be made $H$-free by removing from it at most $\de \binom{n}{2}$ edges.
		\end{theorem}
		
		Now we can deduce Proposition \ref{prop:Erdos-Stone_stability} from Theorems \ref{stability theorem} and \ref{graph removal lemma} as follows.
		
		\begin{proof}[Proof of Proposition \ref{prop:Erdos-Stone_stability}]
			Let $\de=\de_{\ref{stability theorem}}(H,\eps)/2$, $\ga=\min\{\ga_{\ref{graph removal lemma}}(H,\de),\de\}$ and $C=C_{\ref{stability theorem}}(H,\eps)$. Since $G$ contains at most $\ga n^{v(H)}$ copies of $H$, Theorem \ref{graph removal lemma} shows that $G$ contains an $H$-free subgraph $G'$ with ${e(G')\ge e(G)-\de \binom{n}{2}}$. Hence 
			\[
			e(G')\ge \left(\frac{\chi(H)-2}{\chi(H)-1}-\ga-\de\right)\binom{n}{2} \ge \left(\frac{\chi(H)-2}{\chi(H)-1}-\de_{\ref{stability theorem}}(H,\eps)\right)\binom{n}{2}.
			\]
			Moreover, $v(G')=n \ge C = C_{\ref{stability theorem}}(H,\eps)$.
			Therefore, one can apply Theorem \ref{stability theorem}  to obtain a $(\chi(H)-1)$-colourable subgraph $G''$ of $G'$ with
			$v(G'')\ge (1-\eps)n$ and $\de(G'')\ge \left(\frac{\chi(H)-2}{\chi(H)-1}-\eps\right)n$. 
		\end{proof}
		
		\begin{proof}[Proof of Lemma \ref{lem:common_neighbours}]
			Choose $D=qd^{-r}$ and $\rho=e^{-q}d^{rq}$.
			Let $S$ be the set of tuples $(w_1,\ldots,w_r,A)$ where $w_s \in W_{(s)}$ for all $s$, and $A \in \binom{N(w_1,\ldots,w_r)}{q}$.
			We find that
			\begin{equation}
			\label{eq:appendix_1}
			\card{S}=\sum\limits_{A \in \binom{U}{q}}\prod\limits_{s \le r}\card{N(A)\cap W_{(s)}}=\sum\limits_{(w_1,\ldots,w_r)}\binom{\card{N(w_1,\ldots,w_r)}}{q}.
			\end{equation}
			Moreover, our assumption implies that
			\begin{equation}
			\label{eq:appendix_2}
			\sum\limits_{(w_1,\ldots,w_r)}\card{N(w_1,\ldots,w_r)}=\sum\limits_{u \in U} \prod\limits_{s \le r}\deg(u,W_{(s)})\ge d^r\card{U}\cdot \prod\limits_{s \le r}\card{W_{(s)}}.
			\end{equation}
			Note that the function
			\begin{equation*}
			\binom{x}{q}=
			\begin{cases}
			x(x-1)\cdots(x-q+1)/q!   & \text{if $x \ge q-1$}, \\
			0              & \text{if $x < q-1$}.
			\end{cases}
			\end{equation*} 
			is convex. 
			Thus, we can first apply Jensen's inequality to the right hand side of \eqref{eq:appendix_1} and then use the inequality \eqref{eq:appendix_2} to obtain
			$\card{S} \ge \binom{d^r\card{U}}{q} \prod\limits_{s \le r}\card{W_{(s)}}$. We infer from this and the first identity in \eqref{eq:appendix_1} that there is a subset $A \in \binom{U}{q}$ with
			\[
			\prod\limits_{s \le r}\card{N(A)\cap W_{(s)}} \ge \frac{\binom{d^r\card{U}}{q}}{\binom{\card{U}}{q}} \cdot\prod\limits_{s \le r}\card{W_{(s)}} \ge e^{-q}d^{rq} \cdot\prod\limits_{s \le r}\card{W_{(s)}}=\rho \cdot\prod\limits_{s \le r}\card{W_{(s)}},
			\]
			where the second inequality holds since $\binom{\card{U}}{q} \le \left(\frac{e\card{U}}{q}\right)^q$, and 
			$\binom{d^r\card{U}}{q} \ge \left(\frac{d^r\card{U}}{q}\right)^q$ for $\card{U}\ge D=qd^{-r} \ge q$. Hence $\card{N(A)\cap W_{(s)}} \ge \rho\card{W_{(s)}}$ for $s\le r$.
		\end{proof}
		
		\section{Concluding remarks}
		\label{sect:concluding}
		
		Bollob\'as \cite[Corollary 3.5.4]{Bollobas78} showed that every $n$-vertex graph with $\lfloor \frac{n^2}{4} \rfloor + 1$ edges contains cycles of lengths from $3$ up to $\lfloor\frac{n+3}{2}\rfloor$, and thus strengthened the Mantel theorem. Using techniques developed in this paper we can prove the following multipartite version of this result; we omit the details. 
		\begin{theorem}\label{cycle lengths}
			Let $\el \ge 10^{20}$, and let $G=\left(V_1\cup \ldots \cup V_{\el},E\right)$ be a balanced $\el$-partite graph on $n$ vertices such that
			\[
			d(V_i,V_j) \ge \tfrac{1}{2} \quad \text{for $i \ne j$}.
			\] 
			Then, $G$ either contains a cycle of length $h$ for each integer $h$ with $3 \le h \le (\tfrac{1}{2}-\frac{2}{\sqrt{\el}})n$ or is isomorphic to a graph in $\mathcal{G}_{\el}^{3}$.
		\end{theorem}
		
		\noindent The balanced $\el$-partite graph obtained by taking the disjoint union of ${K_{\el}\bigl(\lfloor\frac{n}{2\el}\rfloor-1\bigl)}$ and  $K_{\el}\bigl(\lceil\frac{n}{2\el}\rceil+1\bigl)$ has edge densities between parts strictly greater than $\frac{1}{2}$. However, every cycle of this graph has length at most $\frac{1}{2}n+2\el=(\tfrac{1}{2}+o(1))n$ provided $\el=o(n)$. Therefore, the bound $(\frac{1}{2}-\frac{2}{\sqrt{\el}})n$ in the above result is asymptotically best possible.  
		
		A book in a graph is a collection of triangles sharing a common edge. The size of a book is the number of triangles. Let $b(G)$ be the size of the largest book in a graph $G$. Generalising Mantel's theorem, Erd\H{o}s \cite{Erdos62} showed that every $n$-vertex graph $G$ with $\lfloor \frac{n^2}{4}\rfloor +1$ edges satisfies $b(G) \ge \frac{n}{6}-O(1)$. The optimal bound $b(G) \ge \lfloor\frac{n}{6}\rfloor$ was obtained independently by Edwards in an unpublished manuscript \cite{Edwards77}, and by Khad\v{z}iivanov and Nikiforov in \cite{KN79}. We wonder whether a similar result holds for balanced multipartite graphs.
		
		\begin{conj}
			For every $\eps>0$, there is a constant $C=C(\eps)$ such that the following holds for $\el>C$. Let $G=(V_1\cup \ldots \cup V_{\el},E)$ be a balanced $\el$-partite graph on $n$ vertices such that 
			\[
			d(V_i,V_j)>\tfrac{1}{2} \quad \textrm{for every $i \ne j$}.
			\] 
			Then, $b(G)>\left(\tfrac{1}{6}-\eps\right)n$.
		\end{conj}
		
		\noindent According to Theorem \ref{thm:joints}, the above conjecture is true for $\eps \ge \tfrac{1}{6}-3^{-18}$.
		
		Assume $H$ is not an almost colour-critical graph. Theorem \ref{thm:turan}(1) tells us that $d_{\el}(H)\ge \frac{\chi(H)-2}{\chi(H)-1}+\frac{1}{(\chi(H)-1)^2(\el-1)^2}$ for every $\el \ge v(H)$. Furthermore, this estimate is tight for $H=K_{1,2}$, as shown in Remark \ref{rmk:tighness}. It would be very interesting to have a characterisation of the equality case. 
		
		Bondy, Shen, Thomass\'e and Thomassen \cite{BSTT06} determined the value of $d_{\el}(K_k)$ in the case when $\el=k=3$, while Pfender \cite{Pfender12} obtained result in the case when $\el$ is large enough in terms of $k$. The value of $d_{\el}(K_k)$ is not known in the remaining cases. Nevertheless, when $\el=k \ge 4$, Pfender \cite{Pfender} proposed the following conjecture (see \cite[Section 5]{Nagy11} for more details).
		
		\begin{conj}
			The critical edge density $d_{k}=d_k(K_k)$ satisfies the following recurrence formula:
			\[
			d_{2}=0, \quad d_{k}^2(1-d_{k-1})+d_{k}-1=0 \ \text{for $k \ge 3$}.
			\]
		\end{conj}
		
		Finally, we emphasise that there are other interesting multipartite versions of the Tur\'an theorem. For instance, Bollob\'as, Erd\H{o}s and Szemer\'edi \cite{BES75} introduced the function $\de_r(n)$ which is the smallest integer so that every $r$-partite graph with parts of size $n$ and minimum degree $\de_r(n)+1$ contains a copy of $K_r$. The exact values of $\de_r(n)$ was determined completely by Haxell and Szab\'o \cite{HS06} (for odd $r$), and Szab\'o and Tardos \cite{ST06} (for even $r$) via topological methods.     
		
		\section*{Acknowledgement}
		The authors would like to thank Shagnik Das, Tibor Szab\'o and an anonymous referee for helpful suggestions that improved the presentation of this paper.

\end{document}